\documentclass[11pt, a4paper, reqno]{amsart}
\usepackage[english]{babel}
\usepackage{graphicx}

\usepackage[utf8]{inputenc}
\usepackage[T1]{fontenc}
\usepackage{amssymb,amsmath}

\usepackage{multicol}

\DeclareMathOperator{\R}{\mathbb{R}}
\DeclareMathOperator{\Z}{\mathbb{Z}}
\DeclareMathOperator{\RP}{\mathbb{R}\mathrm{P}}
\renewcommand{\H}{\mathbb{H}}
\newcommand{\PSL}{\operatorname{PSL}}
\newcommand{\SL}{\operatorname{SL}}
\newcommand{\SO}{\operatorname{SO}}
\newcommand{\Sym}{\operatorname{Sym}}
\newcommand{\Core}{\operatorname{Core}}
\newcommand{\Aff}{\operatorname{Aff}}
\newcommand{\GL}{\operatorname{GL}}
\newcommand{\rf}{\mathscr{R}}

\newtheorem{theorem}{Theorem}[section]
\newtheorem{prop}[theorem]{Proposition}
\theoremstyle{definition}\newtheorem{defn}[theorem]{Definition}
\theoremstyle{definition}\newtheorem{q}[theorem]{Question}
\newtheorem{remark}[theorem]{Remark}
\newtheorem{ex}[theorem]{Example}
\newtheorem{lemma}[theorem]{Lemma}
\newtheorem{corollary}[theorem]{Corollary}

\usepackage{pgf,tikz}
\usepackage{overpic}
\usepackage{mathrsfs}
\usetikzlibrary{arrows}

\usepackage{color}
\usepackage[hypertexnames=false]{hyperref}

\title{Constructing proper affine actions via higher strip deformations}
\author{Ne\v{z}a \v{Z}ager Korenjak}

\begin{document}
\definecolor{xdxdff}{rgb}{0.49019607843137253,0.49019607843137253,1.}
\definecolor{uuuuuu}{rgb}{0.26666666666666666,0.26666666666666666,0.26666666666666666}
\definecolor{qqqqff}{rgb}{0.,0.,1.}
\maketitle
\begin{abstract}We introduce higher strip deformations, which give  a way of constructing affine deformations of discrete free groups in the image of the irreducible representation $\PSL_2\R\to\SO(2n,2n-1)$. We  use the Margulis invariant to find a properness condition for these deformations, showing that each convex cocompact free group admits affine deformations acting properly on $\R^{4n-1}$. Using our method we can also construct proper affine actions of virtually free groups, which we demonstrate by finding proper actions of $C_2 \star C_3$ on $\R^7.$\end{abstract}
\section{Introduction}

We will discuss a construction of proper affine actions of free groups on affine spaces. These actions give rise to complete affine manifolds, the classification of which is not yet well-understood.
\subsection{Background}
The still-open Auslander conjecture, posed in 1964, asks if all complete compact affine manifolds have virtually polycyclic fundamental groups. The question is motivated by a famous result by Bieberbach, proving that all fundamental groups of compact Euclidean manifolds are virtually abelian.  In response to the Auslander conjecture Milnor \cite{mil} asked if the assumption of compactness is crucial. More specifically, are there free groups acting affinely on $\R^n$ with (noncompact) manifold quotients? In \cite{M}, Margulis constructed the first examples of free groups acting properly discontinuously on $\R^3.$ The spaces obtained as quotients of free groups acting affinely and properly discontinuously on $\R^3$ are called \emph{Margulis spacetimes}.

Using a completely different approach, in \cite{drumm1} and \cite{drumm}, Drumm and Drumm-Goldman construct further affine actions of finitely generated free groups. Indeed,  any discrete finitely generated free group $\Gamma$ in $\SO(2,1)$ - the group of orientation-preserving linear transformations of $\R^3$ preserving a bilinear form of signature $(2,1)$ -  can be the linear part of a proper affine action on affine three-space. This is shown via directly constructing fundamental domains bounded by crooked planes - piecewise linear surfaces in $\R^3$ with surprising disjointness properties.

Given a finitely generated discrete subgroup $\Gamma < \SO(2,1),$ Danciger-Guéritaud-Kassel in \cite{arc} and \cite{dgk} give a precise description of the cone of $\Gamma$-cocycles $u \colon \Gamma \to \R^3$ determining all its proper affine actions. We will call such cocycles \emph{proper affine deformations} of $\Gamma.$ They are parametrized by the arc complex on the hyperbolic surface $\Gamma \backslash \H^2,$ with each cocycle obtained via a construction called an \emph{infinitesimal strip deformation}.
It is shown that crooked planes, and thus fundamental domains of these affine actions, can be recovered from the data of an infinitesimal strip deformation.

In higher dimensions, however, the situation is less well-understood. It is shown by Abels, Margulis and Soifer in \cite{ams} that  a Zariski dense subgroup of $\SO(p,q)$ admitting proper affine deformations exists exactly when $\{p,q\} = \{2n,2n-1\}$ for some $n$. In \cite{s}, Smilga constructs examples of generalized Schottky free groups in $\SO(2n,2n-1)$ admitting proper affine deformations. If some contracting dynamics conditions on the linear part are satisfied, he shows the existence of an open cone of proper affine actions, utilizing a similar technique as Drumm via constructing fundamental domains for these actions. In a slightly different, but related, context of projective geometry Burelle and Treib in \cite{jp} also introduce higher-dimensional projective crooked planes, bounding fundamental domains in projective space.

Studying Margulis spacetimes via constructing fundamental domains is not the only tool available, and not the one we will use here. Already in \cite{M}, Margulis introduced an invariant now called the \emph{Margulis invariant}. Let $\Gamma$ be a discrete subgroup of $\SO(2,1)$ containing only loxodromic elements.  Let  $u\colon \Gamma \to \R^{2,1}$ be a $\Gamma-$cocycle, determining translation parts of each element; the affine action that $u$ determines on $\R^3$ is $\gamma \cdot x = \gamma x + u(\gamma).$ Then define $\alpha_u \colon \Gamma \to \R$ at $\gamma$ as the inner product of the translation vector $u(\gamma)$ and the \emph{neutral vector,} the (appropriately normed and oriented) eigenvalue 1 eigenvector of $\gamma$.
 Margulis showed the opposite-sign lemma; if $\alpha_u$ takes both positive and negative values, the affine action determined by the cocycle $u$ fails to be proper. The same proof first given for dimension 3 works in exactly the same way for higher dimensions.

Goldman, Labourie, and Margulis in \cite{GLM} show that positivity of Labourie's diffused Margulis invariant \cite{L} or the uniform positivity of the \emph{normed} Margulis invariant also gives a sufficient condition for properness of the affine action. Namely, let $\Gamma < \SO(2,1)$ be a discrete finitely generated group containing only loxodromic elements, and let $l(\gamma)$ denote the translation length of $\gamma \in \Gamma.$ Then, if $\frac{\alpha_u(\gamma)}{l(\gamma)} > c > 0$ for some fixed $c$ and all $\gamma \in \Gamma,$ the affine action determined by $u$ is proper.

The Margulis invariant can be defined in the same way for loxodromic elements in $\SO(2n,2n-1)$. The analogous theorem connecting the Margulis invariant to proper affine actions for groups in $\SO(2n,2n-1)\ltimes \R^{2n,2n-1}$ with Anosov linear part was proved by Ghosh-Treib and Ghosh in \cite{aa} and \cite{g}.

For more background and recent results on the topic of proper affine actions, see the survey by Drumm, Danciger, Goldman, and Smilga, \cite{ddgs}.
\subsection{Summary of results}

Here, we will generalize the infinitesimal strip deformation approach from \cite{arc} to work in higher dimensions. Let $\Gamma < \PSL_2\R$ be a discrete finitely generated free group such that $ S =\Gamma \backslash \H^2$ is a convex cocompact hyperbolic surface. Using the irreducible representation $\sigma_{4n-1} \colon \PSL_2\R \to \SO(2n,2n-1) < \SL_{4n-1}\R \cong \PSL_{4n-1}\R,$ we can study $\sigma_{4n-1}(\Gamma)$ by studying the geometry of the surface $S$. We will define \emph{higher strip deformations} in order to construct proper affine deformations of $\sigma_{4n-1}(\Gamma) < \SO(2n,2n-1)$. We will compute their Margulis invariants, and thus obtain some criteria for when such a deformation is proper.
Using the hyperbolic surface $S$, we show
\begin{theorem}[\ref{th}]Let $S = \Gamma \backslash \H^2$ be a noncompact convex cocompact surface. Then $\sigma_{4n-1}(\Gamma)$ admits an open cone of cocycles determining proper affine actions on $\R^{2n,2n-1}.$\end{theorem}
A \emph{strip system} is the data of a collection of properly embedded arcs $\underline a$ on $S,$   a point $p_i$ on each arc  $a_i \in \underline a$ called the \emph{waist},  angles $\theta_i \in [-\pi, \pi],$ and a real number $w$. We turn this data into an affine deformation called a \emph{higher strip deformation} of $\sigma_{4n-1}(\Gamma)$, by explicitly describing the translation part for each $\gamma \in \Gamma.$

\begin{figure}[h!]
\begin{overpic}[scale=0.27, percent]{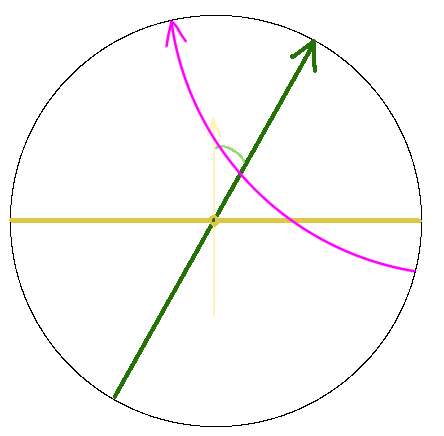}
\put(25,80){$\gamma$}
\put(65,66){$\eta$}
\put(15,50){$a$}
\end{overpic}
\caption{When $\gamma$ crosses the arc $a$, we add the neutral vector of $\sigma_{4n-1}(\eta)$ weighted by $w$ to the translational part of $\sigma_{4n-1}(\gamma).$}\label{f1}
\end{figure}

Here is a brief description of how to build a higher strip deformation along one arc $a$ on $S$. For details, see Section \ref{strips}. Let $c$ be a curve on $S$ representing the element $\gamma \in \Gamma,$ and suppose it crosses $a$. Choose a lift $\tilde a$ of $a$ on $S \cong \H^2,$ and a lift $\tilde c$ of $c$ crossing $\tilde a.$ Denote by $\eta$ the unit-speed hyperbolic translation with axis of translation making the angle $\theta$ through the lift $\tilde p \in \tilde a$ of the waist $ p \in  a$ with the positive perpendicular direction to $\tilde a$, oriented so that $\tilde c$ crosses $\tilde a$ positively. For each such intersection of a curve representing $\gamma$ with the arc $a,$ we add  to the translational part of $\gamma$ the neutral vector of $\sigma_{4n-1}(\eta),$ weighted by the real number $w$. A local picture of a lift about one intersection point is illustrated in Figure \ref{f1}.

For $n = 1$, the accidental isomorphism $\R^{2,1} = \mathfrak{sl}_2\R$ allows us to interpret higher strip deformations as infinitesimal deformations of the hyperbolic surface $S$. If we choose $\theta = 0,$ we recover the infinitesimal strip deformations from \cite{arc}. As the name suggests, these are infinitesimal versions of strip deformations. A strip deformation on a hyperbolic surface $S$ along a properly embedded arc $a$ is a new hyperbolic structure obtained from $S$ by cutting the surface open along $a$ and  gluing in a hyperbolic strip of width $w.$ For $\theta = \pm \frac\pi2$ we get the data of a right (or left) infinitesimal earthquake.

Given a strip system, we can compute  the Margulis invariant of any higher strip deformation obtained from this system for any element $\gamma \in \Gamma$, see Proposition \ref{invt1}. It is a sum of contributions coming from all of the intersections of the arcs with $\gamma$. At each intersection, the contribution depends on the relative position of the axes of $\gamma$ and $\eta$, and can be expressed as a rational function in the cross-ratio of the fixed points in $\partial \H^2$ of $\gamma$ and $\eta.$ For $n = 1,$ the contribution function only has one root; if the angle between $\eta$ and $\gamma$ is less than $\frac\pi2,$ the contribution is positive. In particular, this always happens for $\theta = 0$ and $w > 0,$ independently of the choice of the waist. This shows that the sum of positive orthogonal strips for $n=1$ - the infinitesimal strips from \cite{arc} - determine proper actions if we sum over filling arcs, as in that case all elements of $\Gamma$ have uniformly positive Margulis invariants.

For $n > 1$, however, the contribution function upon one crossing of a given higher strip deformation has multiple simple real roots, changing sign $2n$ times. Even for $\theta = 0$ and $w>0$, we can get both positive and negative contributions, depending on how far from $p$ and at what angle the axis of $\gamma$ crosses $a$.  Consequently, it is not in general the case that any choice of waist would contribute positively to the Margulis invariant. Nonetheless, we show that a careful choice of $p$ and $\theta$ for a convex cocompact surface $S =\Gamma \backslash \H^2$ restricts the possible relative positions of the axes of $\eta$ with the axes of $\gamma \in \Gamma$ enough that we still achieve properness for an open cone of cocycles for any $n$.
In analogue with Smilga's contraction condition on the linear part in \cite{s}, we show that when the convex core is "thin enough", $\theta = 0$ works for any choice of waist.

We can also utilize our approach to construct affine actions of virtually free groups on $\R^{2n,2n-1}.$ We demonstrate this by the example of the free product $C_2 \star C_3$ of cyclic groups of orders $2$ and $3$.  We focus on this example because $C_2 \star C_3$ is abstractly isomorphic to $\PSL_2\Z,$ a group of independent interest. Note that the image of $\PSL_2\Z$ under $\sigma_3$ and $\sigma_7,$ or the image of a convex cocompact deformation of $\PSL_2\Z,$ has one-dimensional first group cohomology, and thus only one affine deformation up to scale and coboundaries. We can realize this cocycle as a higher strip deformation. For $n = 1,2,$ the Margulis invariants of all the loxodromic elements in $\sigma_{4n-1}(\PSL_2\Z)$ are positive. Because the group contains parabolic elements and the results in \cite{GLM} and \cite{aa}, \cite{g} only apply in the purely loxodromic case, this does not a priori mean that the action is proper, though we believe it to be. A generalization of the Goldman-Labourie-Margulis properness criterion to the case containing parabolic elements is under investigation.

Choosing a convex cocompact embedding of $C_2 \star C_3$ in $\PSL_2\R$ - "opening up the cusp" of the modular surface - ensures that all infinite-order elements are loxodromic. Now, our Margulis invariant computation for $n=1,2$ shows that for this deformed copy of $\PSL_2\Z$, the affine deformation we get from a higher strip deformation is proper. In the case of $\R^{2,1}$ and $\R^{4,3}$ this construction describes all possible affine deformations. In higher dimensions, the space of affine deformations is larger, and not all of them come directly from a single strip deformation. However, we can still explore the space of affine deformations using superpositions of higher strips, but we do not do so here.

\subsubsection{Further questions}As mentioned, for $\PSL_2\Z$ acting on $\R^{2n,2n-1},$ the single strip deformation does not account for all possible affine deformations when $n \ge 3$. This leads us to believe that there are some (as of yet unknown) missing parameters that could help in better describing the cone of proper affine cocycles. Some options include assigning translational parts to triangles in an ideal triangulation of $\Gamma \backslash \H^2$ and considering "superpositions" of several higher strip deformations along the same arc.

As another interesting byproduct of the $\PSL_2\Z$ computation, we obtain explicit examples of a free group $\Gamma < \PSL_2\R$ such that $\Gamma \backslash \H^2$ is a three-holed sphere, and an affine deformation in $\R^{6,5}$ such that all three cuff curves have positive Margulis invariant, but the action is not proper, with some non-simple curve achieving negative Margulis invariant.
 This is in stark contrast to the three-dimensional case, where for a three-holed sphere group, requiring positivity of the cuff curves is enough to ensure properness of an affine deformation, see for instance \cite{cdg} and \cite{arc}. For $n>1,$ the space of affine deformations of $\sigma_{4n-1}(\Gamma)$ is more than three-dimensional, so the cone of proper affine deformations, which is always a proper convex cone, cannot be determined by just three inequalities coming from the cuff curves. However, we can still ask
\begin{q}Let $\Gamma \backslash \H^2$ be a three-holed sphere. Is the cone of proper affine deformations of $\sigma_{4n-1}(\Gamma)$ determined by finitely many inequalities?
\end{q}

This example also raises a more general question; in three dimensions, if the Margulis invariant of a cocycle $u \colon \Gamma \to \R^{2,1}$ takes both positive and negative values, $\sup_{\gamma \in \Gamma}\alpha_u(\gamma)$ is approached by a sequence of \emph{simple} curves, see \cite{maxlam}, \cite{glmm}. Because of that, in the three-dimensional case checking only the Margulis invariants of simple closed curves is enough to determine properness. We know that at least for the three-holed sphere, that is no longer sufficient. The question, then, is
\begin{q} Is there a "nice" subset $P_n$ for each $n$ of $\pi_1S$ such that, for any $\pi_1S$-cocycle $u$ the affine action it determines on $\R^{2n,2n-1}$ is guaranteed to be proper as soon as $\frac{\alpha_u(\gamma)}{l(\gamma)}\ge c > 0$ for all $\gamma \in P_n$?\end{q}

\subsection{Organization of the paper} In Section \ref{preliminaries} we recall some facts about hyperbolic and affine geometry, as well as fix 
notation and introduce the version of the Margulis invariant we'll be working with. In Section \ref{strips} we define \emph{higher strip 
deformations}, motivated by strip deformations in \cite{arc}. We compute the Margulis invariants of these building blocks for proper affine 
actions in Proposition \ref{invt1}. In Section \ref{aps} we show that some subset of higher strip deformations always gives us an open cone of 
proper affine actions, which is the statement of Theorem \ref{th}. We also investigate the necessity of some of our assumptions. In the last section, \ref{virt}, we make some observations about actions of $\PSL_2\Z.$

\subsection*{Acknowledgments} The author would like to thank their advisor, Jeff Danciger, for all of his help, encouragement,  insights, and nigh-infinite patience in guiding them through this project. We would also like to thank Jean-Philippe Burelle and Bill Goldman for expressing interest and being willing to talk about some of these results in their initial stages, and François Guéritaud for feedback.

\section{Preliminaries}\label{preliminaries}
In this section, we recall some facts about hyperbolic and affine geometry, as well as fix some notation. 

\subsection{Hyperbolic geometry}\label{sub1}

We denote by $\H^2$ the hyperbolic plane. When we have to make explicit computations, we will use either the Poincaré disk model or the upper half-plane model. In both cases, we think of the boundary $\partial \H^2$ as a copy of $\RP^1$. The group of isometries of $\H^2$ is isomorphic to $\PSL_2\R,$ and its action on the Poincaré disk model or the upper half-plane model is by M\"obius transformations. The action of $\PSL_2\R$ on $\H^2$ extends to the boundary.

We call an element $\gamma$ of $\PSL_2\R$ \emph{loxodromic} if its representative $\tilde \gamma$ in $\SL_2\R$ has two distinct real eigenvalues. In that case, the action of $\gamma$ on $\H^2 \cup \partial \H^2$ has two fixed points in the boundary; equivalence classes in $\RP^1$ of the eigenvectors of $\tilde \gamma.$ We denote the attracting fixed point by $\gamma^+$ and the repelling fixed point by $\gamma^-.$

An element $\gamma$ with a repeated real eigenvalue is \emph{parabolic} and has one fixed point on the boundary, and an element with non-real eigenvalues is \emph{elliptic} with a fixed point in $\H^2.$

The \emph{cross-ratio} of four distinct points $z_1, z_2, z_3, z_4$ in $\partial \H^2 \cong \RP^1$ is the number $$[z_1, z_2; z_3, z_4] := \frac{z_3 - z_1}{z_3 - z_2} \frac{z_4 - z_2}{z_4 - z_1}. $$

The cross-ratio is invariant under the action of $\PSL_2\R,$ meaning that for any $\gamma \in \PSL_2\R,$ we have $[\gamma \cdot z_1, \gamma \cdot z_2 ; \gamma \cdot z_3, \gamma \cdot z_4] = [z_1, z_2; z_3, z_4].$ Further, in our chosen convention, $[t, 1; 0, \infty] = t.$ 

If the pairs $z_1, z_2$ and $z_3, z_4$ are intertwined, meaning that the geodesic in $\partial \H^2$ connecting $z_1$ and $z_2$ intersects the geodesic between $z_3$ and $z_4,$ we can express  $$[z_1,z_2; z_3,z_4] = \frac{\cos(\phi) + 1}{\cos(\phi) - 1},$$ where $\phi$ is the angle between the two aforementioned  geodesics.

An oriented topological surface $S$ of negative Euler characteristic can be endowed with a hyperbolic structure. We do this by realizing $\Gamma = \pi_1(S)$ as a discrete subgroup of $\PSL_2\R$ acting on $\H^2$ by deck transformations so that $\Gamma \backslash \H^2$ is homeomorphic to $S$. The \emph{limit set} $\Lambda_\Gamma$ of $\Gamma$ is $$\Lambda_\Gamma = \overline{\Gamma \cdot x} \cap \partial \H^2 \subseteq \partial \H^2,$$ where $x$ is any point in $\H^2$. It is a closed subset of the boundary, invariant under the action of $\Gamma$ and independent of the choice of $x \in \H^2.$

Denote by $C(\Lambda_\Gamma)$  the convex hull of the limit set of $\Gamma$ inside of $\H^2.$  The set $\Core(S) :=\Gamma \backslash C(\Lambda_\Gamma)$ is called the \emph{convex core} of $S.$ Any geodesic joining two points in $\Core(S)$ lies entirely in $\Core(S).$ In particular, any geodesic representative of a (free) homotopy class of a closed curve lies in $\Core(S).$ If $\Core(S)$ is compact, we call the surface $S$, and the group $\Gamma,$ \emph{convex cocompact.} The condition on the core is equivalent to saying that all elements of $\Gamma \setminus \{1\}$ are loxodromic.

\subsection{Affine geometry}
One of our main objects of interest is the (real) affine space, modeled on $\R^n.$ When there is no chance of confusion, we will abuse notation and denote $n$-dimensional affine space by $\R^n.$ The group preserving the affine structure on $\R^n$ is the group of \emph{affine transformations,} $ \Aff(\R^n) \cong \operatorname{GL}_n\R \ltimes \R^n,$ with $(A, u) \in \operatorname{GL}_n\R \ltimes \R^n$ acting on $x \in \R^n$ by $(A, u)\cdot x = Ax + u.$ We call $A$ the \emph{linear part} and $u$ the \emph{translation}. $\Aff(\R^n)$ comes equipped with two natural projections, $L \colon \Aff(\R^n) \to \operatorname{GL}_n\R,$ which picks out the linear part of an affine transformation, and $u \colon \Aff(\R^n) \to \R^n,$ assigning the translation part to an affine transformation.

For a group $\Gamma < \GL_n\R,$ a \emph{$\Gamma$-cocycle} is a map $u \colon \Gamma \to \R^n$ satisfying $u(\gamma_1 \gamma_2) = \gamma_1\cdot u(\gamma_2) + u(\gamma_1).$ This condition ensures that $\Gamma_u := \{(\gamma, u(\gamma))\} \subset \Aff(\R^n)$ is still a group. Further, the linear part of $\Gamma_u$ is the group $\Gamma$ we started with; $L(\Gamma_u) = \Gamma$. We will often call $\Gamma_u < \Aff(\R^n)$ an \emph{affine deformation} of $\Gamma < \GL_n\R.$

An \emph{affine $n$-manifold} is an $n-$manifold equipped with an affine atlas, meaning that all charts map to the affine space $\R^n,$ and the transition maps are restrictions of affine transformations (on each connected component). A \emph{complete affine manifold} is an affine manifold obtained as a quotient of affine space, $\Gamma \backslash \R^n,$ where $\Gamma < \Aff(\R^n)$ is a discrete subgroup acting properly discontinuously on $\R^n.$

\subsection{Irreducible representations of $\PSL_2\R$} Up to conjugation, there is only one irreducible representation $\tilde\sigma_n \colon \SL_2\R \to \GL_n\R$ for every $n$. One of our main objects of study here are discrete subgroups in the image of the representation $\sigma_n$, so we will need to work with it directly a lot. Let us fix a representative and some notation.

We get the irreducible representation via the action of $\SL_2\R$ on the symmetric power $\Sym^{d-1}(\R^2) \cong \R^d$. We will think of $\Sym^{d-1}(\R^2)$ as a vector subspace of $(\R^2)^{\otimes (d-1)}$ spanned by the symmetric tensors. 

Fix an oriented basis $\{e_1, e_2\}$ of $\R^2,$ and a volume form $\omega$ on $\R^2$ such that $\omega(e_1, e_2) = \omega(e_2, e_1)= 1$ and $\omega (e_i,e_i) = 0.$ As a volume form, $\omega$ is preserved by $\SL_2\R.$ Then the basis for $\Sym^{d-1}(\R^2)$ we'll be working with is
\begin{eqnarray}\label{basis}e_1^{d-1} & := &\underbrace{ e_1 \otimes e_1 \otimes \cdots \otimes e_1}_{d-1} \nonumber \\
e_1^{d-1}e_2 &:=&  e_1 \otimes e_1 \otimes \cdots \otimes e_1 \otimes e_2 +  e_1  \otimes \cdots \otimes e_2 \otimes e_1 + \cdots + e_2 \otimes e_1 \otimes e_1 \cdots \otimes e_1 \nonumber \\
&\vdots& \\
e_1^ie_2^{d-1-i}&:=& \sum_{\substack{(k_1, k_2, \ldots, k_{d-1})\\ \text{ exactly $i$ of the } k_j \\ \text{are $1$, the rest are $2$}}}e_{k_1}\otimes e_{k_2}\otimes \cdots \otimes e_{k_{d-1}} \nonumber \\
e_2^{d-1} &:=& e_2 \otimes e_2 \otimes \cdots \otimes e_2. \nonumber
\end{eqnarray}
Similarly, for arbitrary vectors $v_1, v_2 \in \R^2,$ we denote by $$v_1^iv_2^{n-1-i} := \sum_{\substack{(k_1, k_2, \ldots, k_{d-1})\\ \text{ exactly $i$ of the } k_j \\ \text{are $1$, the rest are $2$}}}v_{k_1}\otimes v_{k_2}\otimes \cdots \otimes v_{k_{d-1}}.$$
For $\gamma \in \SL_2\R,$ we define the map $\tilde \sigma_d \colon \SL_2\R \to \GL_d\R$ by defining it on the basis as $\sigma_d(\gamma)\cdot e_1^ie_2^{d-1-i} = (\gamma e_1)^i (\gamma e_2)^{d-1-i}.$ If $\gamma \in \SL_2\R$ has eigenvalues $\lambda, \lambda^{-1}$ with eigenvectors $v_1, v_2,$ then $\tilde\sigma_d$ has eigenvalues of the form $\lambda^{d-2k-1}$ with corresponding eigenvector $v_1^{d-k-1}v_2^k,$ for each $0 \le k \le d-1.$

For odd $d = 2n-1$, this map descends to a map $\sigma_d \colon \PSL_2\R \to \PSL_d\R.$ Moreover, $\sigma_{2n-1}$ factors through $\SO(n,n-1),$ as each $\sigma_{2n-1}(\gamma)$ preserves a bilinear form $B$:

Define a bilinear map $B \colon (R^2)^{\otimes (d-1)}\times  (R^2)^{\otimes (d-1)} \to \R $ by defining it on elementary tensors as $$B(v_1 \otimes \cdots \otimes v_{d-1}, w_1 \otimes \cdots \otimes w_{d-1}) =- \Pi_{i=1}^{d-1} \omega(v_i, w_i).$$ If we restrict $B$ to $\Sym^{d-1}(\R^2),$ we get a symmetric bilinear form of signature $(n, n-1)$. Because $\SL_2\R$ preserves $\omega,$ the group $\sigma_d(\PSL_2\R)$ preserves $B$ and thus lies in $\SO(B).$

In the basis described in \eqref{basis}, the matrix representing $B$ has $(-1)^k\binom{2n-2}{k-1}$ on the $k$th anti-diagonal entry and zeros everywhere else; 
$$B = \left( \begin{matrix}
0 & 0 & 0 & \cdots & 0 & -1 \\
0& 0 & 0 & \cdots & 2n-2 & 0 \\
\vdots & \vdots & \vdots & \reflectbox{$\ddots$} & \vdots & \vdots \\
0 & 0 & -(n-1)(2n-3) & \cdots & 0 & 0 \\
0 &2n-2 & 0 & \cdots & 0 & 0 \\
-1 & 0 & 0 & \cdots & 0 & 0
 \end{matrix} \right).$$

We will be looking for proper affine deformations of groups $\Gamma < \SO(n,n-1),$ which can only exist if $n$ is even. From now on, we thus restrict to $\SO(2n,2n-1),$ where the signature $(2n,2n-1)$ form we're preserving is the aforementioned $B$. When $\Gamma < \PSL_2\R \to \sigma_{4n-1}(\PSL_2\R),$ we will often not distinguish between $\Gamma$ as a subgroup of $\PSL_2\R$ and $\Gamma$ as a subgroup of $\SO(2n,2n-1)$ in our notation. Similarly, we will confuse $\gamma$ and $\sigma_{4n-1}(\gamma)$ if we believe there is little cause for confusion.

\subsection{The Margulis invariant}The Margulis invariant is a tool for detecting properness of affine actions, first introduced by Margulis in \cite{m2}, \cite{M} in the context of affine actions on affine 3-space. It has been extended to work for groups with linear part in $\sigma_{4n-1}(\PSL_2\R) < \SO(2n,2n-1)$ by Goldman-Labourie-Margulis in \cite{GLM} and for Anosov representations $\Gamma \to \SO(2n,2n-1)$ in \cite{aa}, \cite{g}.

\begin{defn}\label{nv}Let $\gamma \in \PSL_2\R$ be a loxodromic element with eigenvalues $|\lambda| > |\lambda^{-1}|$ and corresponding eigenvectors $v_1, v_2.$ Assume $(v_1, v_2)$ is a positively oriented basis. Then the \emph{neutral vector} of $\sigma_{4n-1}(\gamma)$ is $$x^0(\sigma_{4n-1}(\gamma)) := \frac{1}{B(v_1^{2n-1}v_2^{2n-1},v_1^{2n-1}v_2^{2n-1})^{\frac12}}v_1^{2n-1}v_2^{2n-1}.$$
\end{defn}
\begin{remark}This vector really is well-defined, as we require $v_1$ to correspond to the eigenvalue of larger modulus, and require $(v_1, v_2)$ to be positively oriented. If we choose $-v_1$ instead of $v_1,$ we have to choose $-v_2$ instead of $v_2,$ and the minus signs cancel out to give a vector with the same orientation. Further, we require the vector to have norm $1$, which we can do, as it will always be space-like.

When there's no fear of confusion, we will write $x^0(\gamma)$ instead of $x^0(\sigma_{4n-1}(\gamma)).$
\end{remark}

We are now ready to define the Margulis invariant. Let $\rho \colon \Gamma \hookrightarrow \PSL_2\R \to^{\sigma_{4n-1}} \SO(2n,2n-1)$ be a representation with only loxodromic elements.   Take $u$  an affine deformation of $\Gamma < \SO(2n,2n-1).$ 
\begin{defn}The \emph{Margulis invariant} of $\Gamma_u < \sigma_{4n-1}(\PSL_2\R)\ltimes \R^{2n,2n-1}$ is a function $\alpha_u \colon \Gamma \to \R$ defined by $$\alpha_u(\gamma) = B(x^0(\gamma), u(\gamma)).$$
Here, $x^0(\gamma)$ is the neutral vector of $\gamma$ as in definition \ref{nv}. \end{defn}

We can immediately observe some nice properties of $\alpha_u$:
\begin{itemize}\item $\alpha_{u_1 + u_2} = \alpha_{u_1} + \alpha_{u_2},$
\item $\alpha_{au} = a\alpha_u,$
\item $\alpha_u(\gamma^n) = n\alpha_u(\gamma)$ for $n >0,$
\item $\alpha_u(\eta^{-1}\gamma \eta) = \alpha_u(\gamma)$ for all $\eta, \gamma \in \Gamma,$ that is, $\alpha_u$ is a class function on $\Gamma.$ If $u$ is defined on all of $\PSL_2\R,$ we can take $\eta \in \PSL_2\R.$ 
\end{itemize}

In \cite{M}, \cite{m2} Margulis showed a connection between the Margulis invariant of $\Gamma_u < \SO(2,1)\ltimes \R^{2,1}$ and the properness of the action $\Gamma_u $ on $ \R^{2,1}:$
\begin{theorem}[\cite{M}]Let $\Gamma < \SO(2,1)$ contain only loxodromic elements and let $u$ be a $\Gamma$-cocycle. Suppose there are $\gamma, \eta \in \Gamma\setminus 1$ such that $\alpha_u(\gamma)\alpha_u(\eta) \le 0.$ Then the affine action of $\Gamma_u$ on $\R^{2,1}$ is not proper. 
\end{theorem} 

The theorem works exactly the same for $\Gamma < \sigma_{2n-1}(\PSL_2\R).$ In particular, it shows why we need to restrict to even $n$; if $n$ is odd, we get $\alpha_u(\gamma) = - \alpha_u(\gamma^{-1}),$ so the action cannot be proper in that case.

To get the reverse statement, the correct thing to look at turns out to be the \emph{diffused Margulis invariant,} $\Psi_u,$ first introduced by Labourie in \cite{L}. For $\Gamma < \sigma_{4n-1}(\PSL_2\R)$ it is defined on the space of currents, which are flow-invariant $\Gamma$-invariant probability measures on the unit tangent bundle of the hyperbolic plane, $U\H^2.$ Equivalently, they are flow-invariant measures on $\Gamma \backslash U\H^2$, the unit tangent bundle of the hyperbolic surface $\Gamma \backslash \H^2,$ if we view $\Gamma$ as a subgroup of $\PSL_2\R.$ For our purposes, it is enough to know that free homotopy classes of closed curves embed into the space of currents, and are dense in it. Computing the diffused Margulis invariant on the current $c_\gamma$ representing a closed curve $\gamma$ on $S = \Gamma \backslash \H^2$ gives
$$\Psi_u(c_\gamma) = \frac{1}{l(\gamma)}\alpha_u(\gamma),$$
where $l(\gamma) := \inf_{p \in \H^2}d(p, \gamma \cdot p)$ is the translation length of $\gamma \in \PSL_2\R.$

Using the diffused Margulis invariant, Goldman-Labourie-Margulis in \cite{GLM} slowed that an affine deformation $u$ of $\Gamma < \SO(2n,2n-1)$ is proper if and only if $\Psi_u$ does not take the value $0$. Due to closed curves being dense in the space of currents, and the space of currents being connected, the version of the theorem we will use here is
\begin{theorem}[\cite{GLM}]\label{glm}Let $\Gamma < \PSL_2\R$ be a discrete subgroup consisting only of loxodromic elements. Let $u$ be an affine deformation of $\sigma_{4n-1}(\Gamma).$ Then $\Gamma_u < \SO(2n,2n-1)\ltimes \R^{2n,2n-1}$ acts properly on $\R^{2n,2n-1}$ if and only if there exists $c > 0$ so that $\alpha_u(\gamma) > l(\gamma)\cdot c$ for all $\gamma \in \Gamma$ (or the same is true for -u).
\end{theorem}

\section{Higher strip deformations}\label{strips}
In this section, we define a \emph{higher strip deformation}, a kind of affine deformation of $\Gamma < \sigma_{4n-1}(\PSL_2\R)$ with nice computational properties. The construction relies heavily on the construction of \emph{infinitesimal strip deformations} from \cite{arc}, and recovers it when $n =1.$

Let $S = \Gamma \backslash \H^2$ be a (convex cocompact oriented) hyperbolic surface. We will associate some data to $S$ that we will turn into a $\Gamma$-cocycle $u$.
\begin{defn}A \emph{strip system} on $S$ is a triple $(\underline a, \underline p, \underline \theta),$ where $\underline a$ is a system of disjoint properly embedded transversely oriented geodesic arcs $a_i$ on $S$, $\underline p$ is a set of points $p_i \in a_i,$ and $\underline \theta$ is  a set of angles $\theta_i \in [-\pi, \pi).$\end{defn}
We imagine the angle $\theta_i$ to be based at $p_i,$ with $\theta_i$ being measured from the normal in the positive transverse direction to $a_i$.
\begin{figure}[h!]
\begin{overpic}[scale=0.2, percent]{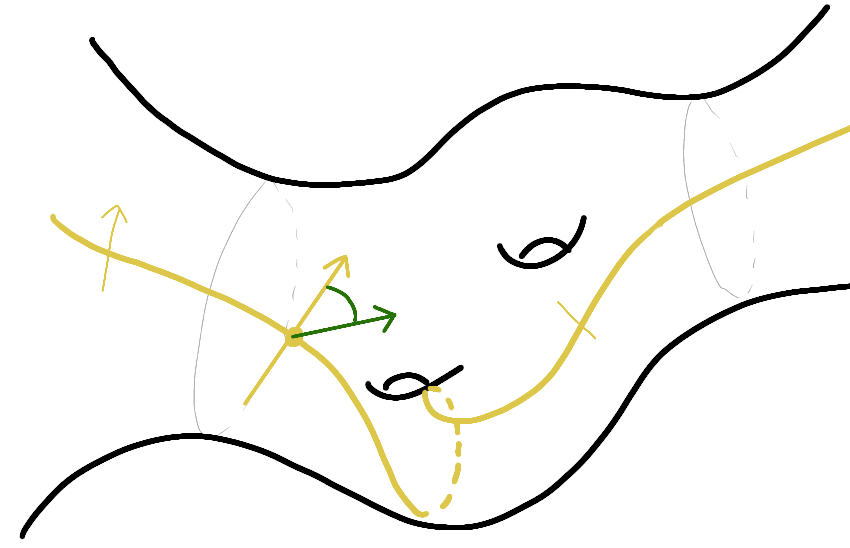}
\put(26,20){$p$}
\put(42,30){$\theta$}
\put(10,40){$a$}
\end{overpic}
\caption{A strip system consisting of one arc on a surface $S$.}
\end{figure}

We now describe how to obtain a cocycle for $\sigma_{4n-1}(\Gamma)$ given a strip system. Denote by $(\tilde{\underline a}, \tilde{\underline p}, \tilde{\underline \theta)}$ the lift of $(\underline a, \underline p, \underline \theta)$ to $\H^2.$ Denote arcs in $\tilde{\underline a}$ by $a_i^j$, with $a_i^j$ being a lift of $a_i,$ and similarly $p_i^j$ being the lift of $p_i$ lying on $a_i^j.$ Set $\eta_i^j \in \PSL_2\R$ to be the loxodromic element with translation axis intersecting $a_i^j$ at $p_i^j$ and making angle $\theta_i$ with $a_i^j.$ For each $i,$ choose a \emph{weight} $r_i \in \R,$ denoting the collection of the weights by $\underline r.$ If $\underline r$ is omitted, assume all the $r_i$ equal $1$.

Choose a basepoint $x_0 \in \H^2$ disjoint from $\tilde{\underline a}.$ To describe the cocycle $u_{(\underline a, \underline p, \underline \theta)},$ we describe what its value is for each $\gamma \in \Gamma.$ Let $c_\gamma\colon [0,1]\to \H^2$ be a path connecting $c_\gamma(0)=x_0$ and $c_\gamma(1) =\gamma \cdot x_0$ in $\H^2,$ intersecting each arc in $\tilde{\underline a}$ transversely. For each intersection $y \in \tilde{\underline a} \cap c_\gamma$, let $i(y)$ and $j(y)$ be the indices $i, j$ such that $y \in a_i^j.$ Denote by $\iota_y(\tilde{\underline a}, c_\gamma)$ the intersection number at $y$ of $\tilde{\underline a}$ and $c_\gamma,$ which is $1$ if the orientation of  $c_\gamma$ at $y$ agrees with the transverse orientation of $ a^{j(y)}_{i(y)},$ and $-1$ if the orientations disagree. 

\begin{defn}\label{strip}The map $u_{(\underline a, \underline p, \underline \theta, \underline r)} \colon \Gamma \to \R^{2n,2n-1}$ defined by $$u_{(\underline a, \underline p, \underline \theta, \underline r)}(\gamma) = \sum_{y \in \tilde{\underline a} \cap c_\gamma}r_{i(y)}\iota_y(a_{i(y)}^{j(y)}, c_\gamma)\cdot x^0(\sigma_{4n-1}(\eta_{i(y)}^{j(y)}))$$ is called a \emph{higher strip deformation}.\end{defn}
\begin{figure}[h!]
\begin{overpic}[scale=0.25, percent]{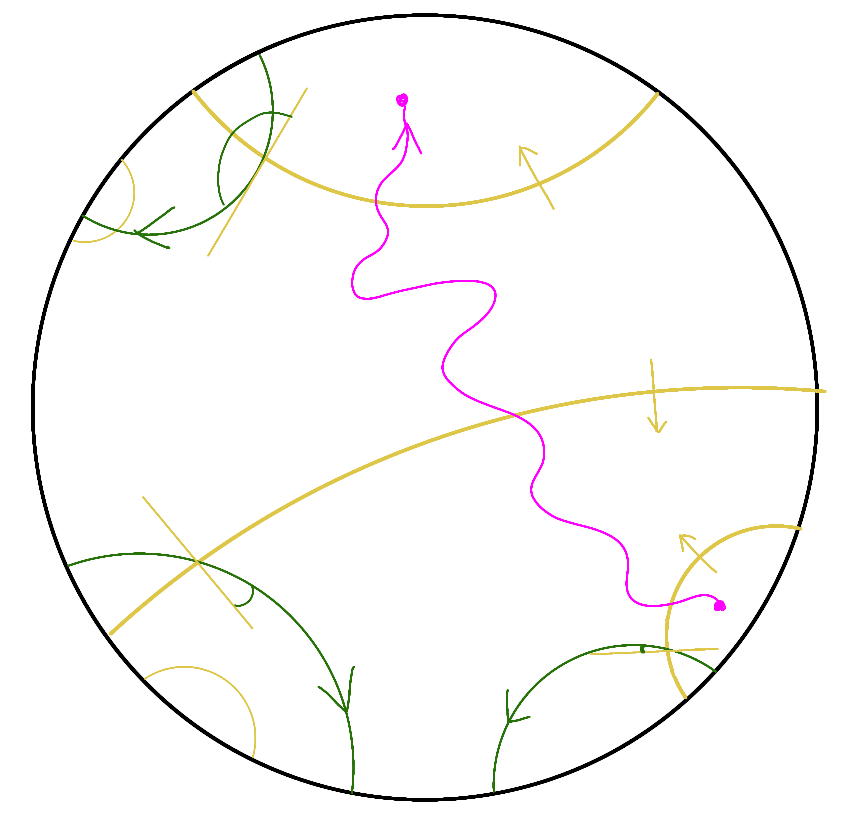}
\put(86,22){\small$x_0$}
\put(49,81){\small$\gamma \cdot x_0$}
\put(88,34){\small$a_1^1$}
\put(60,17){\small$\eta_1^1$}
\put(88,50){\small$a_2^1$}
\put(40,18){\small$\eta_2^1$}
\put(70,70){\small$a_1^2$}
\put(20,61){\small$\eta_1^2$}
\put(60,55){\small$c_\gamma$}
\end{overpic}
\caption{Here, $c_\gamma$ crosses $\tilde{\underline a}$ three times. The cocycle associated to $\gamma$ in this picture would be $u(\gamma) = x^0(\sigma_{4n-1}(\eta_1^1)) - x^0(\sigma_{4n-1}(\eta_2^1)) + x^0(\sigma_{4n-1}(\eta_1^2)).$}\label{coc}
\end{figure}
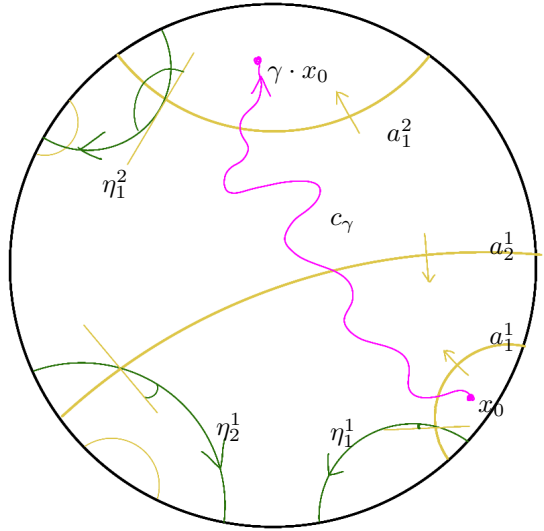
This definition does not depend on the choice of the path $c_\gamma.$ It does depend on the choice of basepoint $x_0,$ but only up to conjugation. Because we will be interested in computing the Margulis invariant of $u_{\underline a, \underline p, \underline \theta},$ and the Margulis invariant is conjugation-invariant, that will not matter for our computations in the future and we can always move the basepoint wherever we want to for any given computation. We can always choose $c_\gamma$ to be a lift of the geodesic representative of $\gamma$'s free homotopy class on $S$.

\begin{remark}The map $u_{(\underline a, \underline p, \underline \theta, \underline r)} \colon \Gamma \to \R^{2n,2n-1}$ defined in Definition \ref{strip} is a cocycle and determines an affine deformation of $\sigma_{4n-1}(\Gamma).$
\proof Due to linearity it is enough to check that a strip system consisting of one arc $a$ with weight $1$ determines a $\Gamma$-cocycle. Let $u = u_{a, p, \theta, 1}$ and choose $x_0 \in \H^2 \setminus \tilde a.$ Let $\gamma_1, \gamma_2$ be two elements of $\Gamma$. We want to show that $u(\gamma_1 \gamma_2) = \sigma_{4n-1}(\gamma_1)u(\gamma_2) + u(\gamma_1).$ Let $c_{\gamma_1\gamma_2}, \, c_{\gamma_1},$ and $c_{\gamma_2}$ be the paths connecting $x_0$ to $\gamma_1 \gamma_2 \cdot x_0,$ $\gamma_1 \cdot x_0,$ and $\gamma_2 \cdot x_0$ respectively. We can choose $c_{\gamma_1\gamma_2}$ to pass through $\gamma_1 \cdot x_0,$ with $c_{\gamma_1\gamma_2(\frac12)}=\gamma_1 \cdot x_0.$ Fixing this choice gives us 
$$u(\gamma_1\gamma_2) = u(\gamma_1) + \sum_{y\in \tilde a \cap c_{\gamma_1\gamma_2}((\frac12,1))}\iota_y(a^{j(y)},c_{\gamma_1\gamma_2})x^0({\sigma_{4n-1}(\eta^{j(y)})}).$$
Because the lift of a strip system is equivariant with respect to $\Gamma,$ we have
$$\iota_y(a^{j(y)}, c_{\gamma_1\gamma_2}) = \iota_{\gamma_1^{-1}\cdot y}(a^{j(\gamma_1^{-1}\cdot y)},\gamma_1^{-1}\cdot c_{\gamma_1\gamma_2}),$$
and $\gamma_1^{-1}\cdot c_{\gamma_1\gamma_2}$ is a path between $x_0$ and $\gamma_2\cdot x_0,$ so (up to reparametrization) we can denote it by $c_{\gamma_2}.$ Again thanks to equivariance, $\eta^{j(y)} = \gamma_1 \eta^{j(\gamma_1^{-1}\cdot y)\gamma_1^{-1}}.$ Each $y' = \gamma_1^{-1}\cdot y$ is a point in the intersection $\tilde a \cap c_{\gamma_2}.$ Thus, the sum $S =\sum_{y\in \tilde a \cap c_{\gamma_1\gamma_2}((\frac12,1))}\iota_y(a^{j(y)},c_{\gamma_1\gamma_2})x^0({\sigma_{4n-1}(\eta^{j(y)})})$ equals

\begin{eqnarray*}S &=& \sum_{y' \in \tilde a \cap c_{\gamma_2}}\iota_{y'}(a^{j(y')}, c_{\gamma_2})x^0(\sigma_{4n-1}(\gamma_1 \eta^{(j(y')}\gamma_1^{-1}))\\
&=&  \sum_{y' \in \tilde a \cap c_{\gamma_2}}\iota_{y'}(a^{j(y')}, c_{\gamma_2})x^0(\sigma_{4n-1}(\gamma_1)\sigma_{4n-1}( \eta^{(j(y')})\sigma_{4n-1}(\gamma_1^{-1}))) =\\
 &=&\sum_{y' \in \tilde a \cap c_{\gamma_2}}\iota_{y'}(a^{j(y')}, c_{\gamma_2})\sigma_{4n-1}(\gamma_1)(x^0(\sigma_{4n-1}( \eta^{(j(y')})) = \\
&=&\sigma_{4n-1}(\gamma_1) \sum_{y' \in \tilde a \cap c_{\gamma_2}}\iota_{y'}(a^{j(y')}, c_{\gamma_2})x^0(\sigma_{4n-1}( \eta^{(j(y')}))= \\ &=& \sigma_{4n-1}(\gamma_1) u(\gamma_2),\end{eqnarray*}
which is what we wanted to show.
\endproof
\end{remark}

\begin{remark}Here, we defined a strip system to consist of arcs, points on arcs, and angles. We could have chosen the strip system to consist of arcs and an equivariant assignment of  elements of $\PSL_2\R$ to each lift of the arcs, avoiding the need for the points and angles, but introducing more annoyingness in keeping track of the lifts. In figure \ref{coc}, that would mean that the green arcs denoting the axes of $\eta_i^j$ need not cross $a_i^j.$ All computations below work just as well in this case.
\end{remark}
\begin{remark}For $n=1$, we recover infinitesimal deformations of surfaces. If we choose all the angles $\theta_i$ to be $0$ or $\pi$, we get an infinitesimal earthquake. If we choose all $\theta_i$ to be $ \frac{\pi}{2},$ we get (positive) infinitesimal strip deformations of speed 1 with waists $p_i$ from \cite{arc}.
\end{remark}
\begin{remark}We note that it is not necessary to require a strip system to consist of disjoint arcs; all definitions would work just as well were the arcs to intersect. We choose to work in the setting where they do not in analogy with infinitesimal strip deformations from \cite{arc}, where there is a non-infinitesimal version of strip deformations, for which disjointness of the arcs is necessary. In future work, we hope to construct a good notion of non-infinitesimal strip deformations in pseudo-Riemannian spaces $\H^{2n,2n-1},$ for which we expect disjointness to be required. 
\end{remark}
The Margulis invariants of infinitesimal strip deformations are always uniformly positive, as shown in \cite{arc}, but the behavior of the Margulis invariant of a higher strip deformation is more mysterious. In the rest of this section, we will explore the basic properties of the Margulis invariant of a higher strip.

Recall that to compute the Margulis invariant of $\gamma,$ we need to compute $u(\gamma)$ and $x^0(\gamma).$ For a higher strip deformation, $u(\gamma)$ is a sum of vectors of the form $x^0(\eta).$ Therefore, it is enough to compute $B(x^0(\gamma), x^0(\eta)),$ which represents the contribution to the Margulis invariant of $\gamma$ when $\gamma$ crosses an arc with associated translation $\eta.$

\begin{prop}\label{invt1} Let $\gamma, \eta \in \PSL_2\R$ be two loxodromic elements. Assume they have distinct fixed points on $\partial \H^2$ and let $t = [\eta^+, \eta^-; \gamma^+, \gamma^-]$ be the the cross-ratio of their attracting and repelling fixed points. Then $$B(x^0(\sigma_{4n-1}(\eta)), x^0(\sigma_{4n-1}(\gamma)) = \frac{1}{(t-1)^{2n-1}}\sum_{j=0}^{2n-1}\binom{2n-1}{j}^2t^j =: \mathscr{R}_n(t).$$
\end{prop}
\proof $B$ is symmetric and invariant under multiplication by elements of $\SO(2n,2n-1).$ Therefore, up to switching the order of $\gamma$ and $\eta$ and up to conjugation, we can assume that $\gamma^+ = t, \, \gamma^- = 1, \, \eta^+ = \infty$ and $\eta^- = 0.$ Then, the attracting eigenvector of $\eta$ is $e_1$ and the repelling eigenvector is $e_2.$ That means that \begin{eqnarray*}x^0(\sigma_{4n-1}(\eta)) &=& \frac{1}{B(e_1^{2n-1}e_2^{2n-1},e_1^{2n-1}e_2^{2n-1})^{\frac12}}e_1^{2n-1}e_2^{2n-1}\\& =& \binom{4n-2}{2n-1}^{-\frac12}e_1^{2n-1}e_2^{2n-1}.\end{eqnarray*} The linear functional $B(e_1^{2n-1}e_2^{2n-1}, -)$ is the projection to the middle coordinate in the basis \eqref{basis}, weighted by  $B(e_1^{2n-1}e_2^{2n-1},e_1^{2n-1}e_2^{2n-1}) = \binom{4n-2}{2n-1}.$ Therefore, in order to compute $B(e_1^{2n-1}e_2^{2n-1}, x^0(\sigma_{4n-1}(\gamma)))$ we only really need the middle coordinate of $x^0(\sigma_{4n-1}(\gamma)).$

 The attracting eigenvector for $\gamma$ is $te_1 + e_2,$ and the repelling eigenvector is $e_1 + e_2.$ If $t > 1,$ this eigenbasis is positively oriented, and if $t< 1,$ it is negatively oriented.

Consider first the case when $t > 1.$ Then as per definition \ref{nv},an appropriately oriented (but not yet normed) eigenvalue $1$ eigenvector for $\sigma_{4n-1}(\gamma)$ is $\tilde x^0(\gamma) = (te_1 + e_2)^{2n-1}(e_1 + e_2)^{2n-1},$ the symmetrization of \begin{equation}\label{neut}\underbrace{(te_1 + e_2)\otimes (te_1 + e_2) \otimes \ldots \otimes (te_1 + e_2)}_{2n-1}\otimes \underbrace{(e_1 + e_2)\otimes \ldots \otimes (e_1 + e_2)}_{2n-1}.\end{equation} 

Write $\tilde x^0(\gamma) = \sum_{k=0}^{4n-2}c_ke_1^{k}e_2^{4n-2-k}.$ We are looking for $c_{2n-1}.$ For the moment, let's consider $\tilde x^0(\gamma)$ as a vector in $(\R^2)^{\otimes 4n-2}$ and write it as 
$$\tilde x^0(\gamma) = \sum b_{i_1, i_2, \ldots, i_{4n-1}}e_{i_1}\otimes e_{i_1}\otimes \cdots \otimes e_{i_{4n-2}}.$$
We have the equality $$c_{2n-1} = b_{\underbrace{1, \ldots, 1}_{2n-1}, \underbrace{2, \ldots, 2}_{2n-1}},$$ so it is enough to compute $b_{1, \ldots, 1, 2, \ldots, 2},$ the coefficient in front of $e_1 \otimes e_1 \cdots e_1 \otimes e_2 \otimes \cdots \otimes e_2.$ It is going to be a polynomial in $t$ of degree at most $2n-1$, so let us write $c_{2n-1} = \sum_{j = 0}^{2n-1}a_jt^j.$

The coefficient $a_j$ is counting the number of tensors in the symmetrization of \eqref{neut} where $j$ of the first $2n-1$ factors are of the form $te_1 + e_2.$ There are $\binom{2n-1}{j}\binom{2n-1}{2n-1-j}$ of those; we choose $j$  among the first $2n-1$ spots to place $te_1 + e_2$ in, and $2n-1-j$ among the latter half of the spots to place the remaining $2n-1-j$ vectors $te_1 + e_2$ in. Therefore, $a_j =\binom{2n-1}{j}\binom{2n-1}{2n-1-j}= \binom{2n-1}{j}^2.$ That gives
$$c_{2n-1} = \sum_{j=0}^{2n-1}\binom{2n-1}{j}^2t^j.$$

We still need to compute the norm of $\tilde x^0(\gamma)$ before we can rest. We compute using the standard volume form $\omega$ on $\R^2$. The only pairs of elementary tensors $te_1 + e_2$ and $e_1 + e_2$ that do not evaluate to $0$ are the ones with $te_1+ e_2$ and $e_1 + e_2$ in exactly the opposite spots. As we have to choose half of the $4n-2$ spots for $te_1 + e_2,$ which determines the rest, there are $\binom{4n-2}{2n-1}$ such pairs, giving us
\begin{eqnarray*}\|\tilde x^0(\gamma)\| &=&\left( -\binom{4n-2}{2n-1}\omega(te_1 + e_2, e_1 + e_2)^{2n-1}\omega(e_1 + e_2, te_1 + e_2)\right)^{\frac12}  = \\
&=& \binom{4n-2}{2n-1}^{\frac12}(-(t-1)^{2n-1}(1-t)^{2n-1})^{\frac12} = \\
&=& \binom{4n-2}{2n-1}^{\frac12}(t-1)^{2n-1}.
\end{eqnarray*}
Combining the computations gives us \begin{eqnarray*}B(x^0(\eta), x^0(\gamma)) &=& B\left(\binom{4n-2}{2n-1}^{-\frac12}e_1^{2n-1}e_2^{2n-1}, \frac{\binom{4n-2}{2n-1}^{-\frac12}}{(t-1)^{2n-1}}\tilde x^0(\gamma)\right)= \\
&=& \binom{4n-2}{2n-1}\binom{4n-2}{2n-1}^{-1} \sum_{j=0}^{2n-1}\binom{2n-1}{j}^2t^j =\\
&=& \frac{1}{(t-1)^{2n-1}}\sum_{j=0}^{2n-1}\binom{2n-1}{j}^2t^j.\end{eqnarray*}

In the case when $t < 1,$ the basis $(te_1 + e_2, e_1 + e_2)$ is negatively oriented, so we have to take $-te_1 - e_2$ instead of $te_1 + e_2.$ Doing the same computation then gives $c_{2n-1} = -\sum_{j=0}^{2n-1}\binom{2n-1}{j}^2t^j.$ Similarly we get $\|x^0(\gamma)\| =  \binom{4n-2}{2n-1}^{\frac12}(1-t)^{2n-1}.$ In the final expression, the minuses cancel and we get the same result. \endproof

We can similarly compute what happens when the fixed points of $\gamma$ and the fixed points of $\eta$ are not all distinct; in that case, we get $B(x^0(\gamma), x^0(\eta)) = 1,$ which is also what happens in the formula of proposition \ref{invt1} as $t$ tends to $\pm\infty$.

\begin{remark}\label{pproperties}The polynomial $\sum_{j=0}^{2n-1}\binom{2n-1}{j}^2t^j$ always has a zero at $-1$. Furthermore, all its zeros are negative, and if $t$ is a zero, so is $\frac1t$. It is a special example of a hypergeometric function, and all its zeros are real and simple, see for instance \cite{poly}. \end{remark}
\begin{lemma}\label{bigt} If the axes of $\gamma$ and $\eta$ do not intersect and have the same orientation, $B(x^0(\sigma_{4n-1}(\eta)), x^0(\sigma_{4n-1}(\gamma)) > 1.$ \end{lemma}
\proof The axes of $\gamma$ and $\eta$ not intersecting and having the same orientation means that the cross-ratio $t = [\eta^+, \eta^-; \gamma^+, \gamma^-]$ is greater than $1$. 
On $t > 1,$ the function $\frac{1}{(t-1)^{2n-1}}\sum_{j=0}^{2n-1}\binom{2n-1}{j}^2t^j$ is strictly decreasing. At $1,$ there is a pole, with $\lim_{t \downarrow 1}\rf_n(t) = \infty.$ As $\lim_{t\to\infty}\rf_n(t) = 1,$ the function $\rf_n$ is strictly greater than $1$ on the interval $(1, \infty).$\endproof

Here are the graphs of $ \frac{1}{(t-1)^{2n-1}}\sum_{j=0}^{2n-1}\binom{2n-1}{j}^2t^j$ for $n = 1, 2, 3.$ Note that for $n=1$ there is only one zero; the sign of the Margulis invariant of an infinitesimal strip deformation depends only on whether or not the transverse orientation of the arc and the orientation of the translation axis of $\gamma$ agree. We expand on this observation in the next section and relate it to \cite{arc}. 
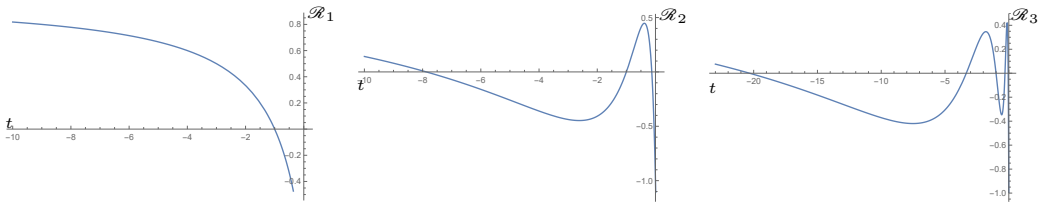
\begin{figure}[h!]
    \centering
    \begin{minipage}{0.3\textwidth}
        \centering
        \begin{overpic}[width=0.9\textwidth]{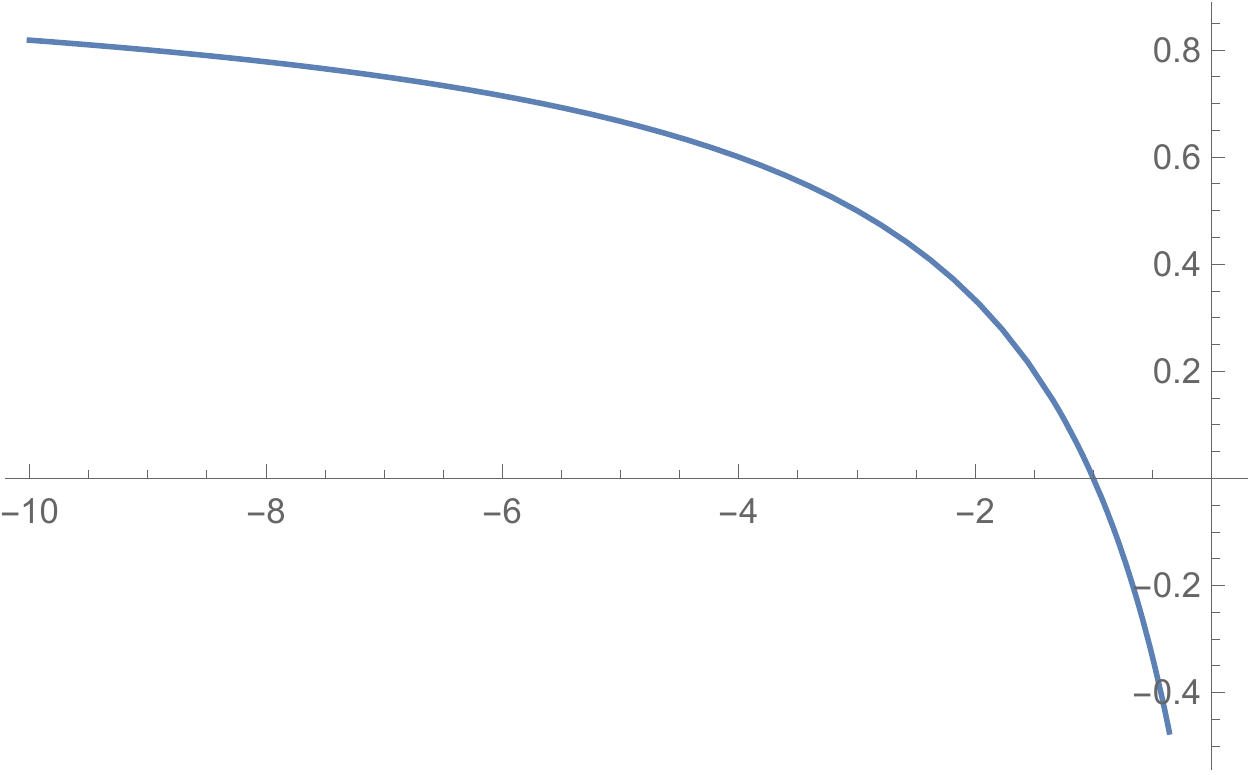}
\put(0,24){\Tiny $t$}
\put(98, 60){\Tiny $\rf_1$}
\end{overpic}
    \end{minipage}
    \begin{minipage}{0.3\textwidth}
        \centering
        \begin{overpic}[width=0.9\textwidth]{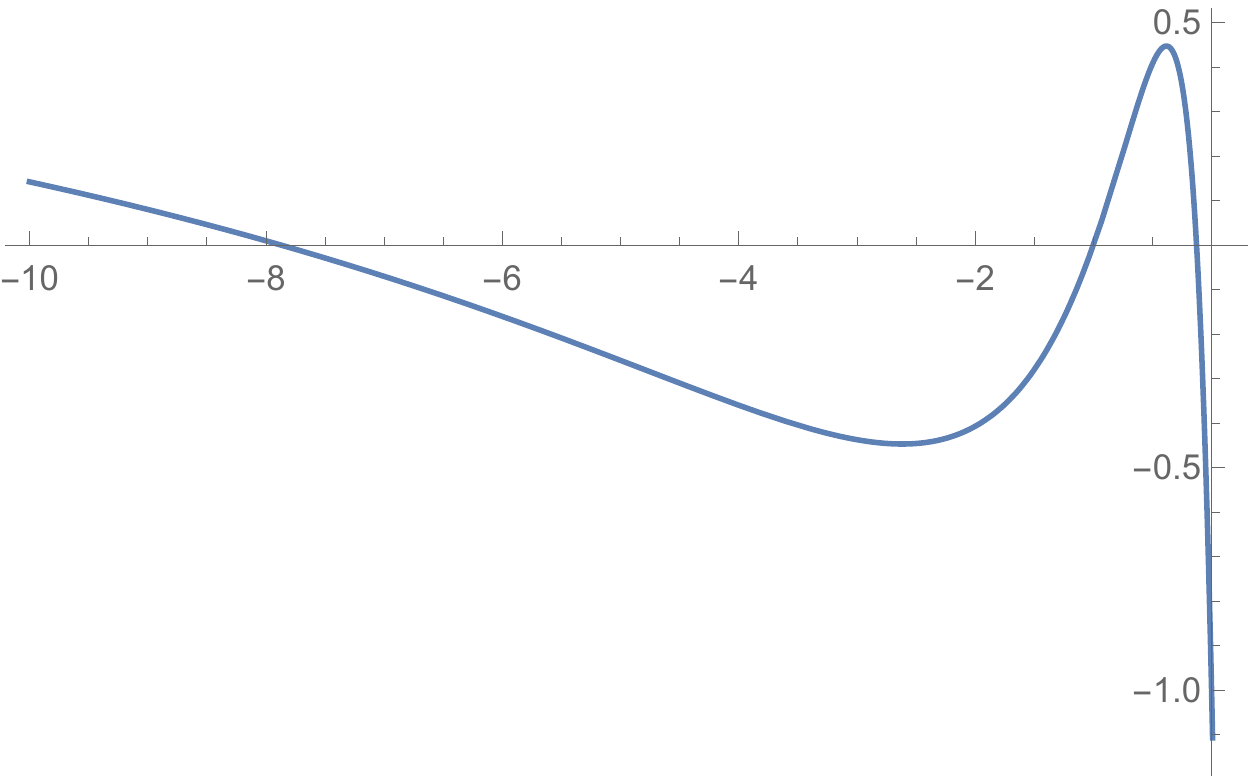} 
\put(98, 60){\Tiny$\rf_2$}
\put(0,37){\Tiny$t$}
\end{overpic}
    \end{minipage}
\begin{minipage}{0.3\textwidth}
        \centering
        \begin{overpic}[width=0.9\textwidth]{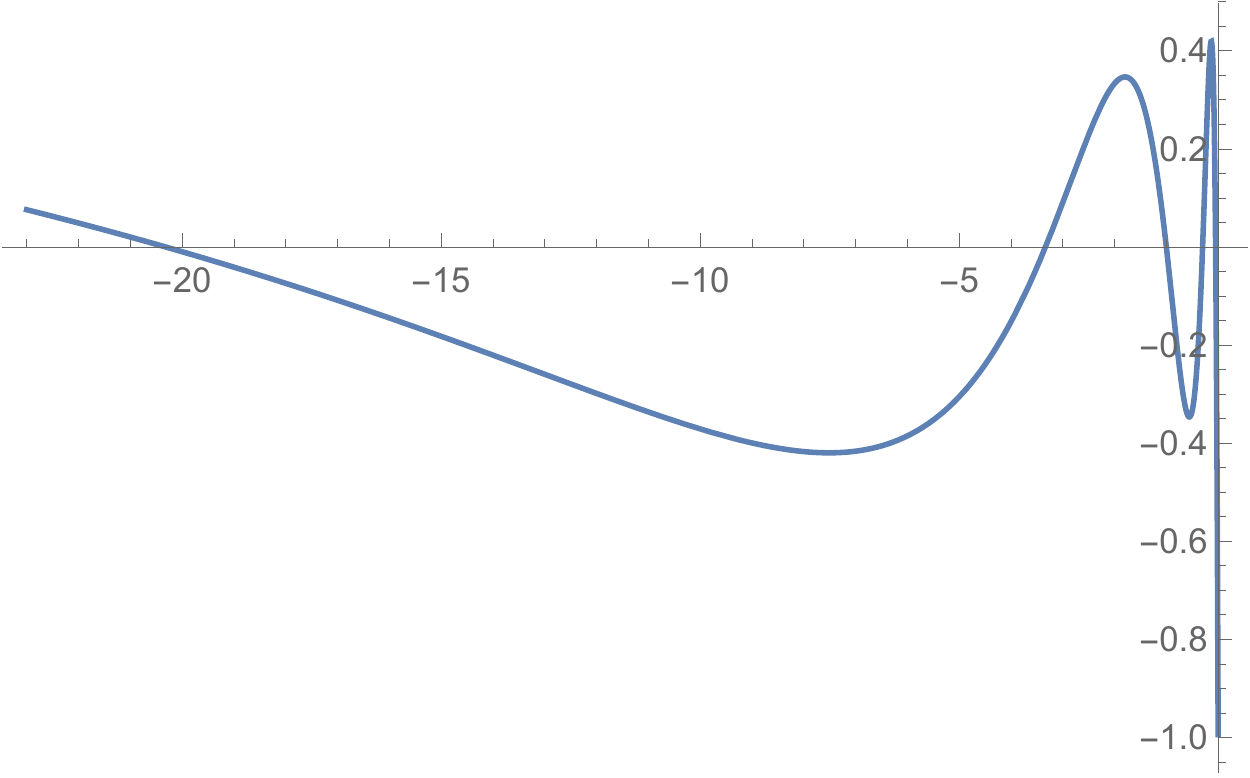} 
\put(0,36){\Tiny$t$}
\put(98, 60){\Tiny$\rf_3$}
\end{overpic}
\end{minipage}
\caption{The Margulis invariant contribution upon one crossing for different $n$, with the $x$-axis being the cross-ratio of the endpoints of $\gamma$ and $\eta$. On the left hand side, the graph for $n=1$ shows that the only zero is at $-1,$ which is where the orientations of $\gamma$ and $\eta$ flip from agreeing to being opposite. In the middle, we have the case for $n=2$, where in addition to $-1,$ two other zeros exist. The rightmost case of $n=3$ shows $5$ zeros. An illustration of the case when $n=2$ from a different point of view is provided in Figures \ref{scha} and \ref{ang}.}
\end{figure}

In higher dimensions the behavior of the Margulis invariant is more complicated; even when orientations agree, we can get either a positive or a negative contribution to the Margulis invariant. Controlling these contributions is thus more difficult and depends more on the geometry of the hyperbolic surface $S$ than it does in the one-dimensional case.

\begin{remark}\label{withangles}When $t \in [-\infty, 0],$ we can express the cross-ratio in terms of the angle between the axes of $\eta$ and $\gamma$, obtaining the notion of "good" and "bad" angles.  For $n = 2,$ we can express the function as $ \frac18(3\cos(\phi) + 5 \cos(3\phi)),$ where $\phi$ is the angle between the axes of $\eta$ and $\gamma$.  The zeros on the interval $(0, \pi)$ are $\frac\pi2$ and approximately $0.684719$ and $2.45687.$ 

We can compute the exact zeros of $\rf_2(t)$ and can explicitly observe the sign-switching behavior in higher dimensions in the following pictures.\end{remark}
\begin{figure}[h!]
    \centering
        \begin{overpic}[scale=0.15, percent]{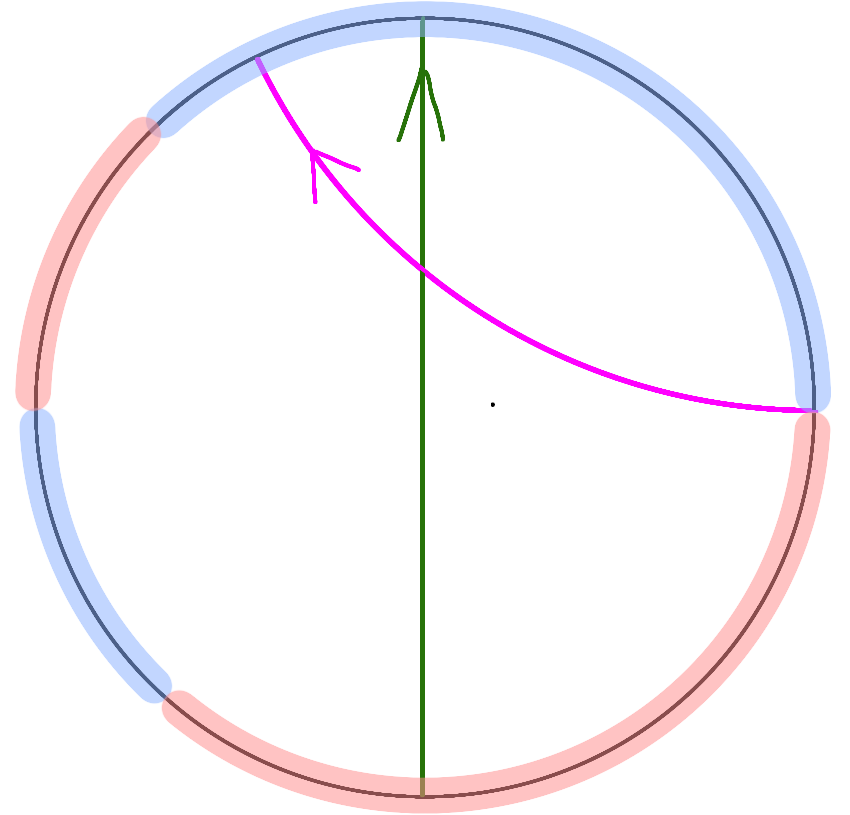} 
\put(42,40){\small$\eta$}
\put(70,54){\small$\gamma$}
\put(46, 97){\small$\infty$}
\put(-9,45){\small$-1$}
\put(48,-6){\small$0$}
\put(100,45){\small$1$}
\put(-17,83){\small$-4-\sqrt{15}$}
\put(-16,5){\small$-4 + \sqrt{15}$}
\put(27,92){\small $t$}
\end{overpic}
   \caption{The green curve $\eta$ is the translation associated to an arc $\gamma$ crosses. For $\R^{4,3},$ the positive contributions to the Margulis invariant of the pink curve $\gamma$ happen if the axis of $\gamma$ points at a blue region, else the contribution is negative, or exactly $0$ at $-1, -4-\sqrt{15},$ and $-4+\sqrt{15}.$}\label{scha}
\end{figure}

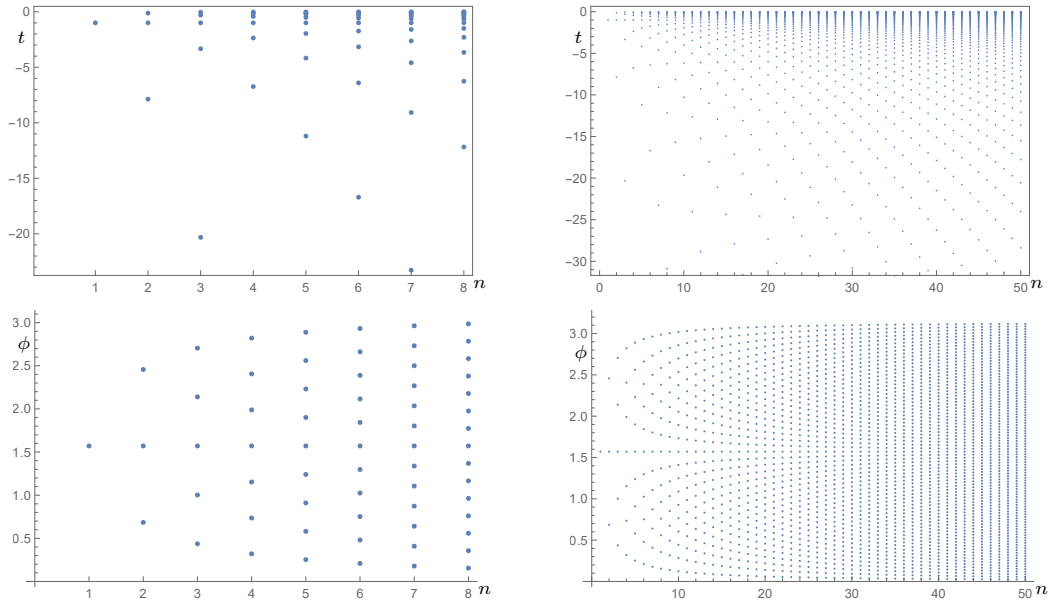
\begin{figure}[h!]
\centering
\begin{minipage}{0.48\textwidth}
\begin{overpic}[width = 0.85\textwidth]{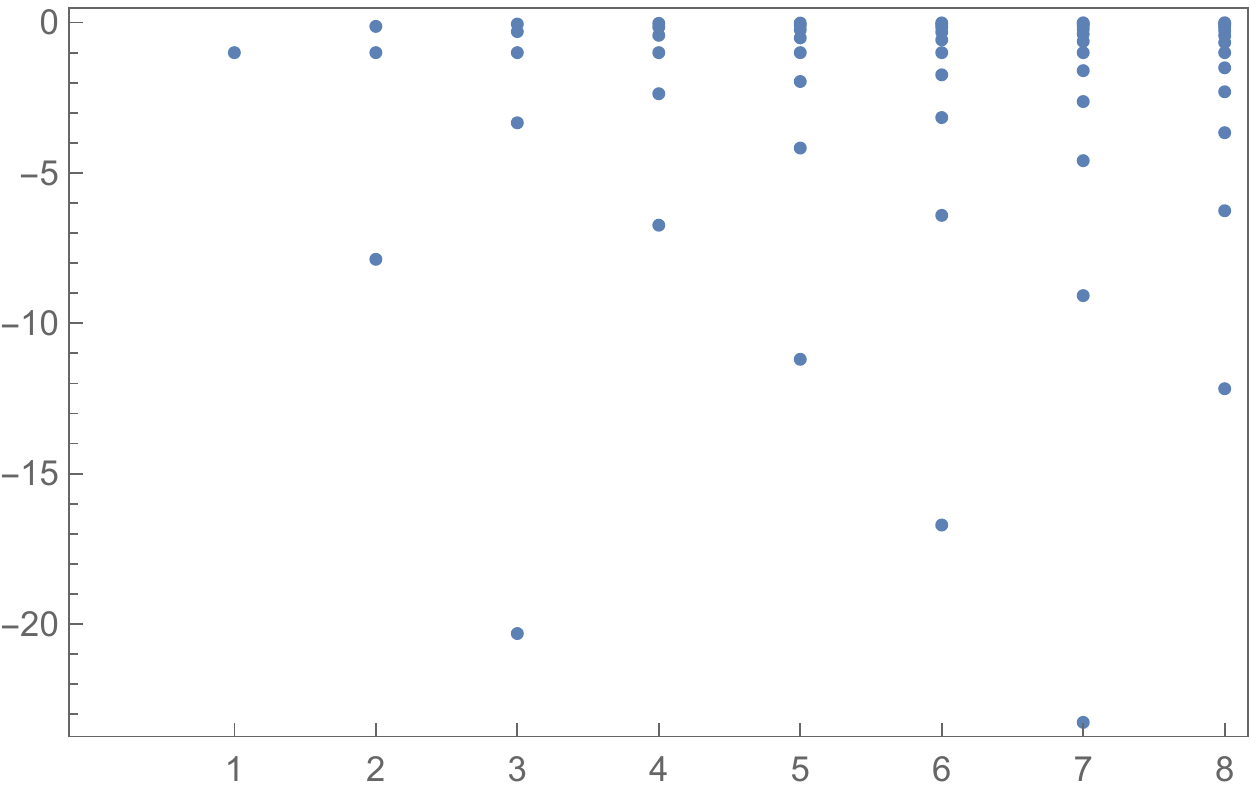}
\put(100,2){\tiny $n$}
\put(2,55){\tiny $t$}
\end{overpic}
\end{minipage}
\begin{minipage}{0.48\textwidth}
\begin{overpic}[width = 0.85\textwidth]{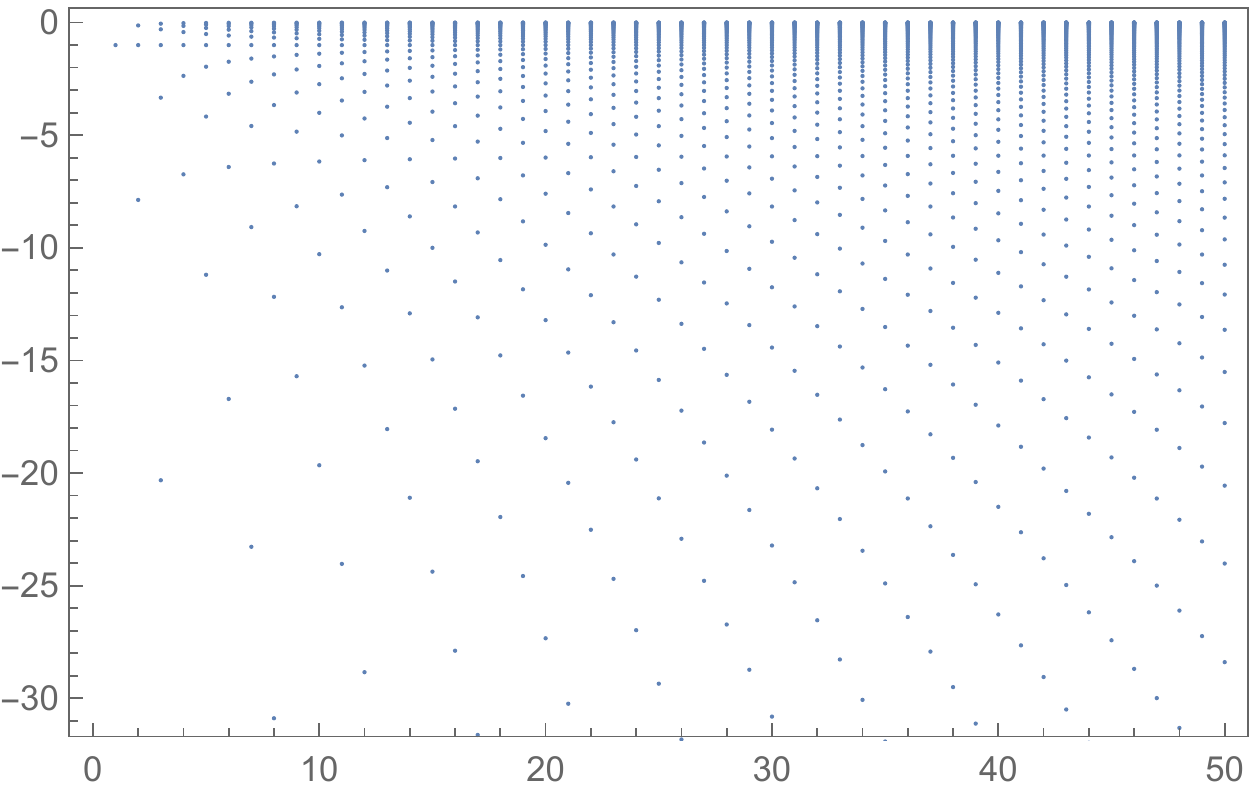}
\put(100,2){\tiny $n$}
\put(2,55){\tiny $t$}
\end{overpic}
\end{minipage}
\vspace{1ex}

\begin{minipage}{0.48\textwidth}
\begin{overpic}[width = 0.85\textwidth]{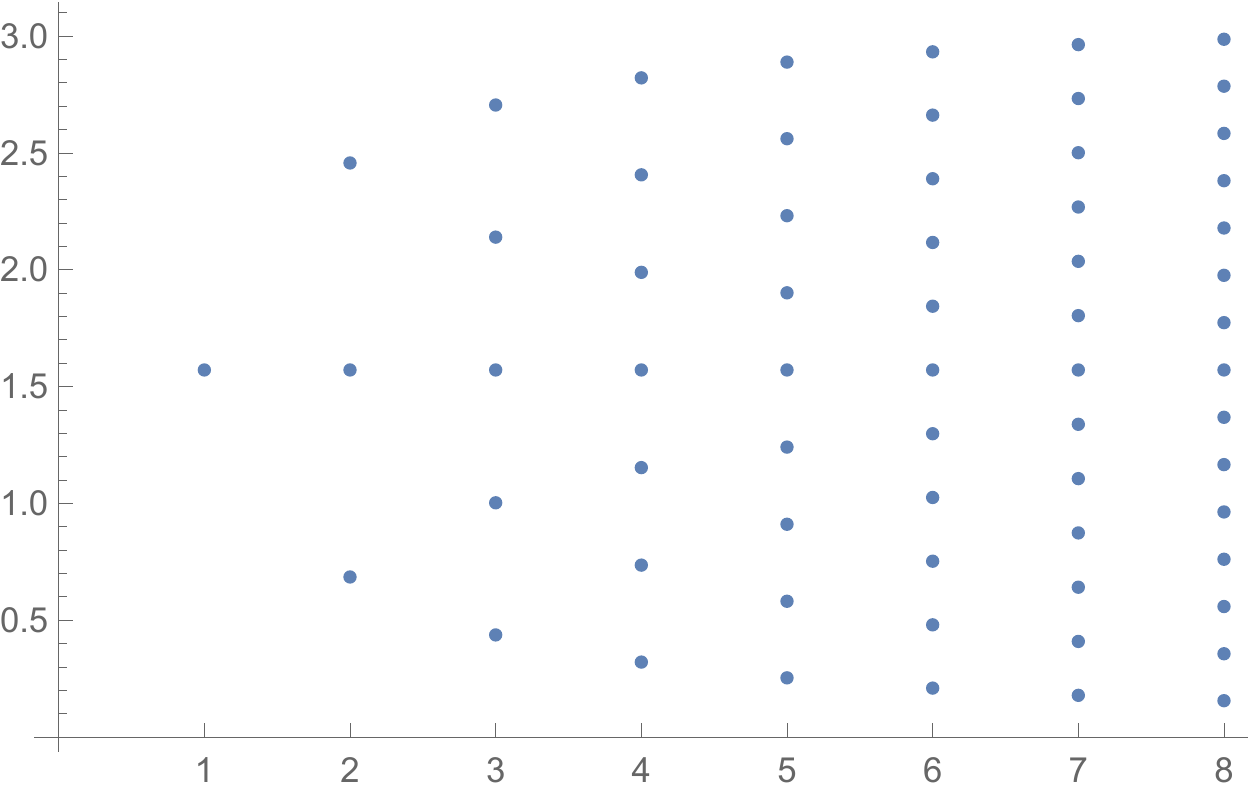}
\put(100,2){\tiny $n$}
\put(1,55){\tiny $\phi$}
\end{overpic}
\end{minipage}
\begin{minipage}{0.48\textwidth}
\begin{overpic}[width = 0.85\textwidth]{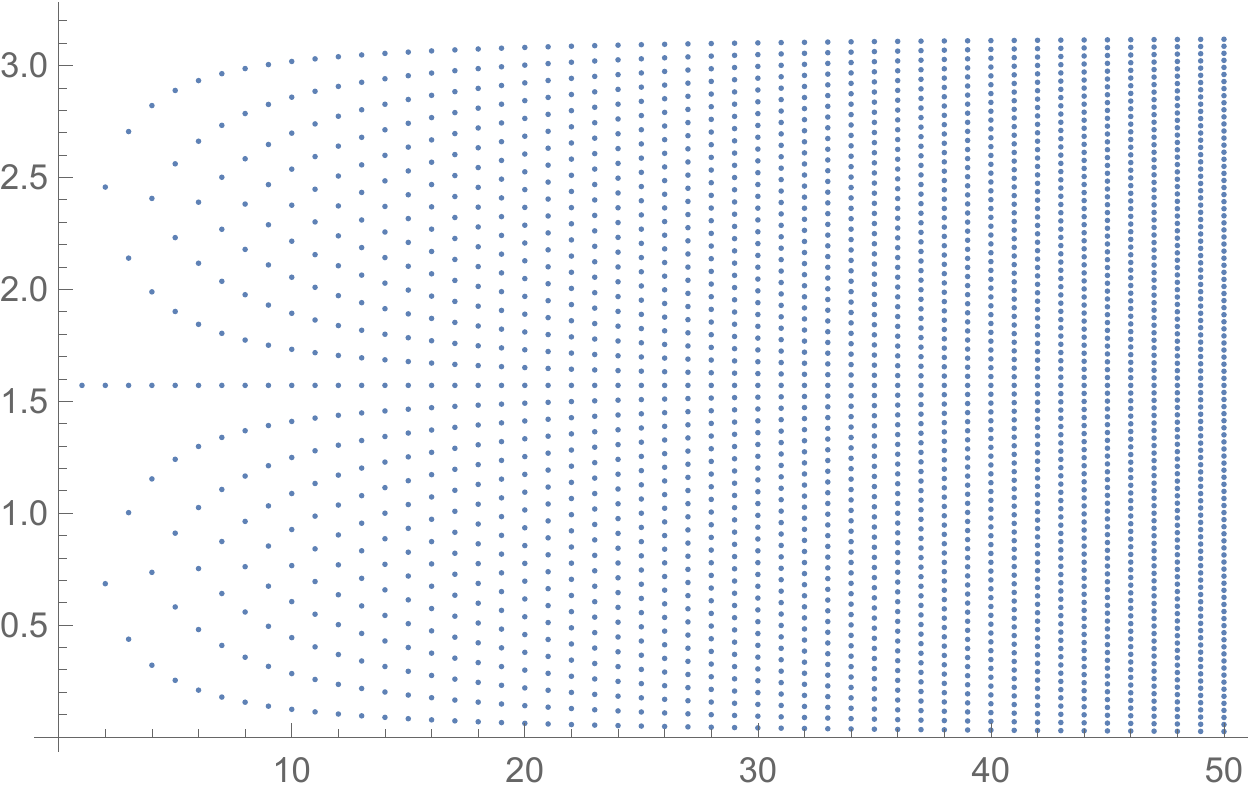}
\put(100,2){\tiny $n$}
\put(1,53){\tiny $\phi$}
\end{overpic}
\end{minipage}
\caption{In the top row, first a graph of the zeros of $\rf_n$ for $n$ between 1 and 8, and on the right, zeros of $\rf_n$ up to $n=50$. In the bottom row, zeros of $\rf_n$ on the interval $(-\infty, 0)$ expressed in terms of the angle $\phi \in (0, \pi)$ between $\gamma$ and $\eta$. On the left zeros up to $n=8$ are pictured, and on the right, $n$ goes up to $50.$}\label{ang}
\end{figure}
\section{Applications}\label{aps}
\begin{figure}[h!]
    \centering
    \begin{minipage}{0.3\textwidth}
        \centering
        \begin{overpic}[width=0.9\textwidth, percent]{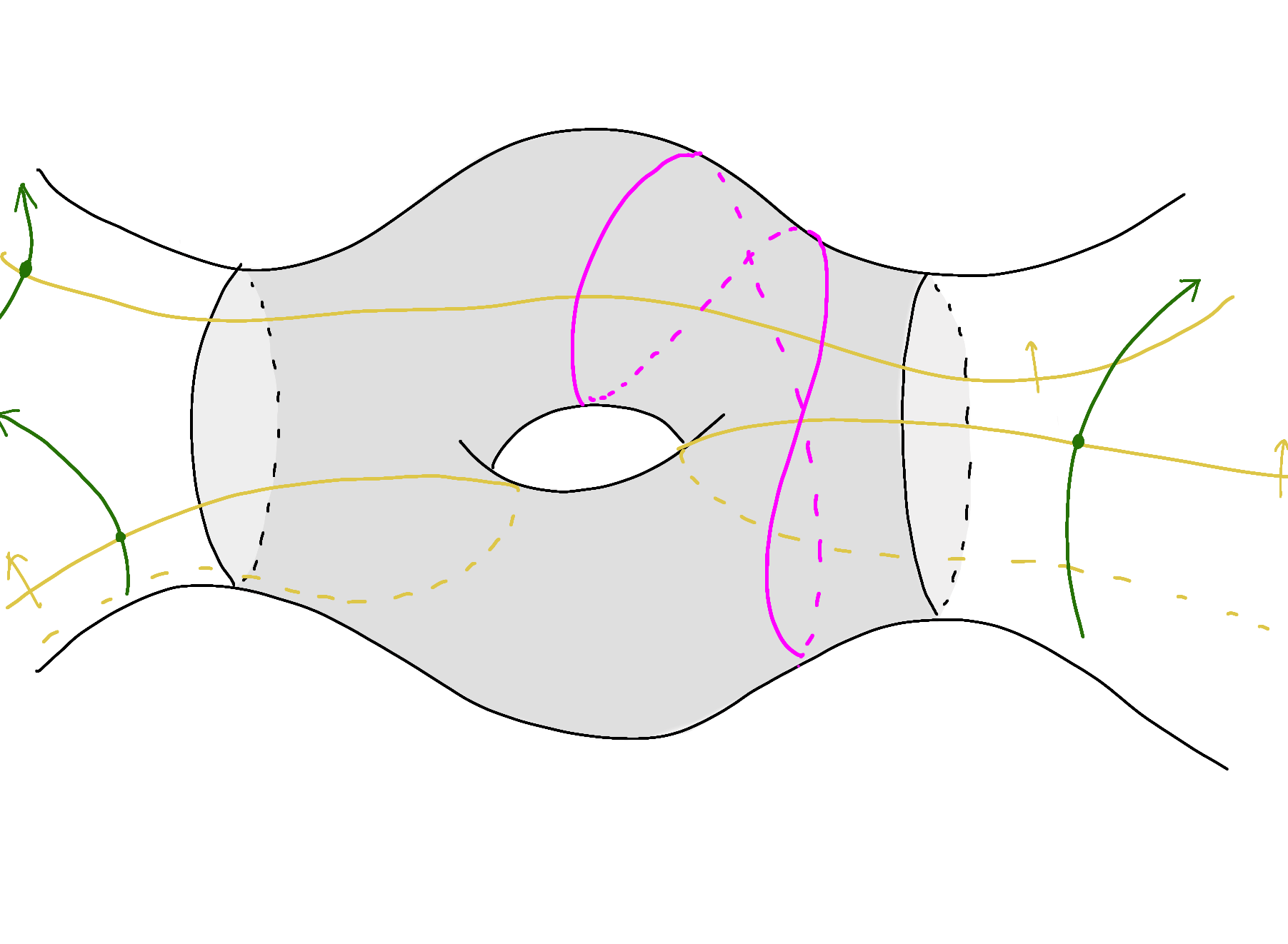}
\put(0,30){\tiny$a_1$}
\put(5,50){\tiny$a_2$}
\put(90,40){\tiny$a_3$}
\put(55,23){\tiny$\gamma$} \end{overpic}
    \end{minipage}
    \begin{minipage}{0.3\textwidth}
        \centering
        \begin{overpic}[width=0.9\textwidth,percent]{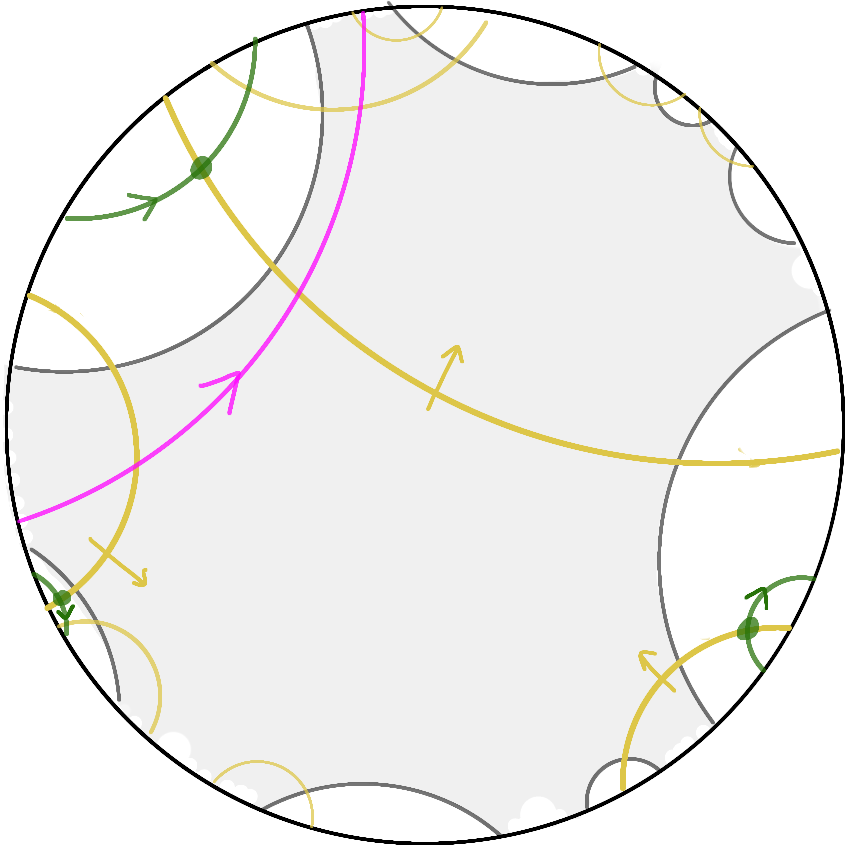} 
\put(28,46){\small$\gamma$}
\put(60,51){\small$a_i^j$}
\put(19,70){\small$\eta_i^j$}\end{overpic}
    \end{minipage}
\begin{minipage}{0.3\textwidth}
        \centering
        \begin{overpic}[width=0.9\textwidth,percent]{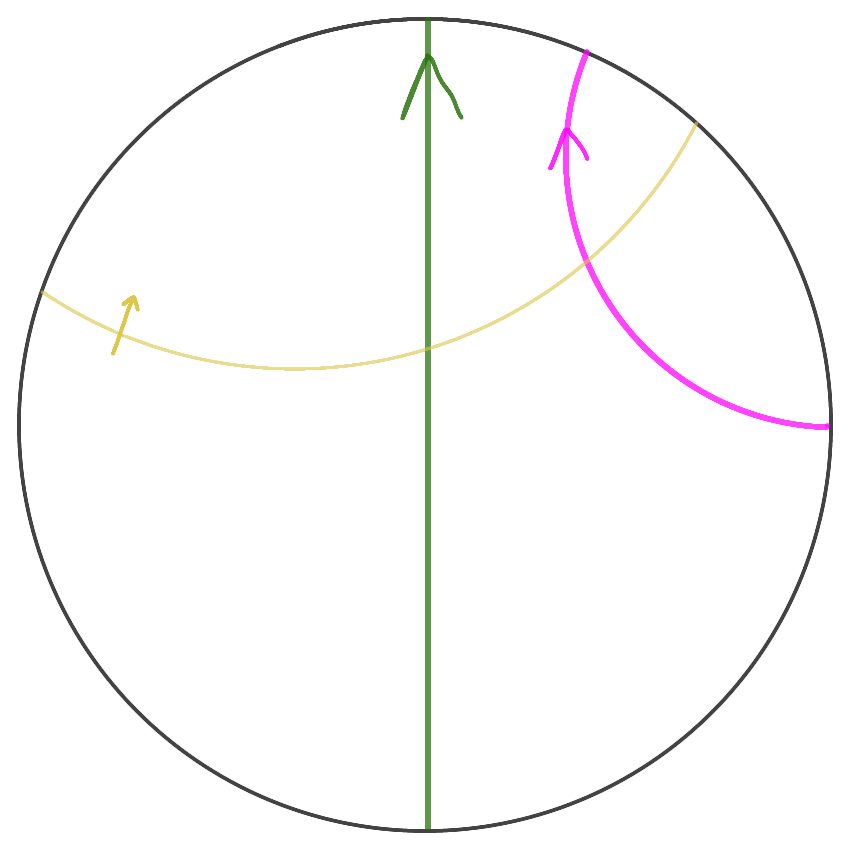}
\put(42,35){\small$\eta_i^j$}
\put(20,50){\small$a_i^j$} 
\put(87,52){\small$\gamma$}
\put(70,95){\small $t$}\end{overpic}
\end{minipage}
\caption{In the leftmost picture, we have a strip system on $S$. As demonstrated in the middle picture, it lifts to the universal cover. In the rightmost picture, we only pay attention to one crossing of an arc in $\tilde{\underline a}$ by the axis of $\gamma$ and the relative position of the translation axes of an $\eta_i^j$ and $\gamma$.}\label{strips2}
\end{figure}
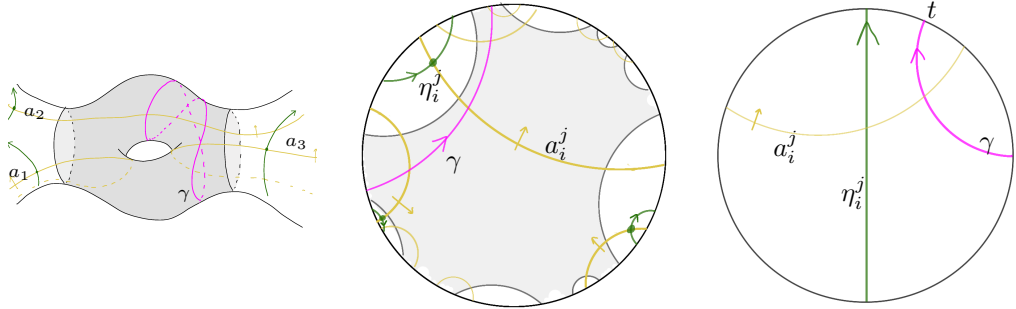

In this section, we use higher strip deformations to show that for a convex cocompact free group $\Gamma < \PSL_2\R,$ the group $\sigma_{4n-1}(\Gamma) < \SO(2n,2n-1)$ always admits proper affine deformations. We also construct examples to compare and contrast the behavior of higher strip deformations to infinitesimal strip deformations, and connect it with some previous work in the field. 
\begin{theorem}\label{th}Let $S = \Gamma \backslash \H^2$ be a convex cocompact noncompact hyperbolic surface. Then there exist proper affine deformations of $\sigma_{4n-1}(\Gamma) < \SO(2n,2n-1).$
\end{theorem}
\proof We will construct a family of proper actions using higher strip deformations. Take a strip system on $S$ such that $\underline a$ is a filling system of arcs, all the points in $\underline p$ lie outside the convex core of $S,$ the angles $\underline \theta$ are chosen so that the geodesic axes of $\eta_i$ through each $p_i$ at transverse angle $\theta_i$ do not intersect the convex core, and all the weights are $1$. Further, let $\theta_ i \in (-\frac\pi2, \frac\pi2),$ so that the orientations of $a_i$ and $\eta_i$ match. A cartoon is provided in Figure \ref{strips2}.

Now let $\gamma \in \Gamma$ be an element. We want to use Theorem \ref{glm} from \cite{GLM}; we need to show that the Margulis invariant $B(x^0(\gamma), u_{(\underline a, \underline p, \underline \theta)}(\gamma))$ is uniformly positive. 

The geodesic representative of $\gamma$ in $S$ lies entirely in the convex core of $S$, and therefore the axis of $\gamma$ in $\H^2$ lies entirely in the lift of the convex core, as demonstrated with the first two pictures in Figure \ref{strips2}. The axes of all the translations $\eta_i^j$ lie entirely outside the convex core. That and the orientations of $a_i^j$ and $\eta_i^j$ agreeing guarantees that $[(\eta_i^j)^+, (\eta_i^j)^-; \gamma^+, \gamma^-] \in (1, \infty].$ A typical relative position of $\gamma$ and $\eta_i^j$ is demonstrated in the third picture of Figure \ref{strips2}. As per Proposition \ref{invt1} and Lemma \ref{bigt}, the contribution to the Margulis invariant of $\gamma$ is at least $1$ every time $\gamma$ crosses an arc in $\underline a.$

Because $S$ is a convex cocompact surface, there is some $0 < R < \infty$ that is the maximum diameter of a piece in $\Core(S)\setminus \underline a$. Because we chose a filling arc system $\underline a,$ every closed geodesic on $S$ will cross an arc in $\underline a$ at least once, and at least once per length $R$. Therefore, we get the contribution of at least $1$ to the Margulis invariant per every $R$ length of $\gamma.$ We can estimate
$$\alpha_{u_{(\underline a, \underline p, \underline \theta)}}(\gamma) \ge \#\text{of crossings of }\underline a\text{ and the geodesic representative of }\gamma \ge \frac{l(\gamma)}{R}.$$
Thus the normed Margulis invariant of $u_{(\underline a, \underline p, \underline \theta)}$ is uniformly positive, and by \ref{glm}, $u_{(\underline a, \underline p, \underline \theta)}$ is a proper affine deformation of $\sigma_{4n-1}(\Gamma).$ \endproof

The proof above is very constructive, so let us restate it in a formulation reflecting that.
\begin{prop}\label{good} Let $S = \Gamma \backslash \H^2$ be a noncompact convex cocompact surface. Let $\underline a$ be a filling arc system on $S$, with $p_i \in a_i$ and $\theta_i$ such that the geodesics through $p_i$ at transverse angles $\theta_i$ do not intersect the convex core of $S$. Let $\underline r$ be a tuple of positive real numbers. Then the strip system $(\underline a, \underline p, \underline \theta, \underline r)$ determines a proper affine cocycle for $\sigma_{4n-1}(\Gamma)$ for every $n$.
\end{prop}

We can use our construction to show that in fact there is an open cone of proper affine deformations:
\begin{theorem}\label{open}Let $S = \Gamma \backslash \H^2$ be a convex cocompact noncompact hyperbolic surface. Then there exists an open cone of proper affine deformations of $\sigma_{4n-1}(\Gamma) < \SO(2n,2n-1)$ spanned by higher strip deformations.
\end{theorem}

\proof Choose generators $\gamma_i$ of $\Gamma$ and an arc system $a_i \in \underline a$ in such a way that each generator $\gamma_i$ intersects only one arc $a_i$ exactly once. An example of such a configuration is provided in Figure \ref{torus}. 

For each $i$, we can independently choose $p_i(t_i) \in a_i$  from an open subinterval of $a_i$ paramterized by $t_i$ - anywhere outside of the convex core - and for each such $p_i(t_i),$ we can choose an open subinterval of angles $\theta_i(t_i, s_i) \subset (\frac\pi2, \frac\pi2)$ so that the associated geodesic does not intersect $\Core(S).$  The one-arc strip systems $(a_i, p_i(t_i), \theta_i(t_i, s_i))_{t,s}$ determine a subset of translations associated to $\sigma_{4n-1}(\gamma_i)$ in $\R^{2n,2n-1}.$ Because $\sigma_{4n-1}$ is an irreducible representation, this subset does not lie in any proper vector subspace of $\R^{2n,2n-1},$ and therefore its positive linear span contains an open cone.
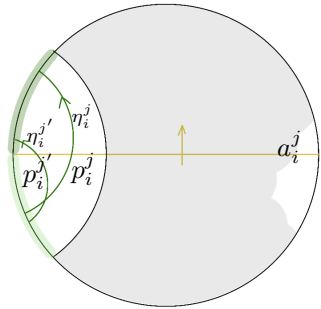
\begin{figure}[h]
\centering
\begin{overpic}[scale=0.22, percent]{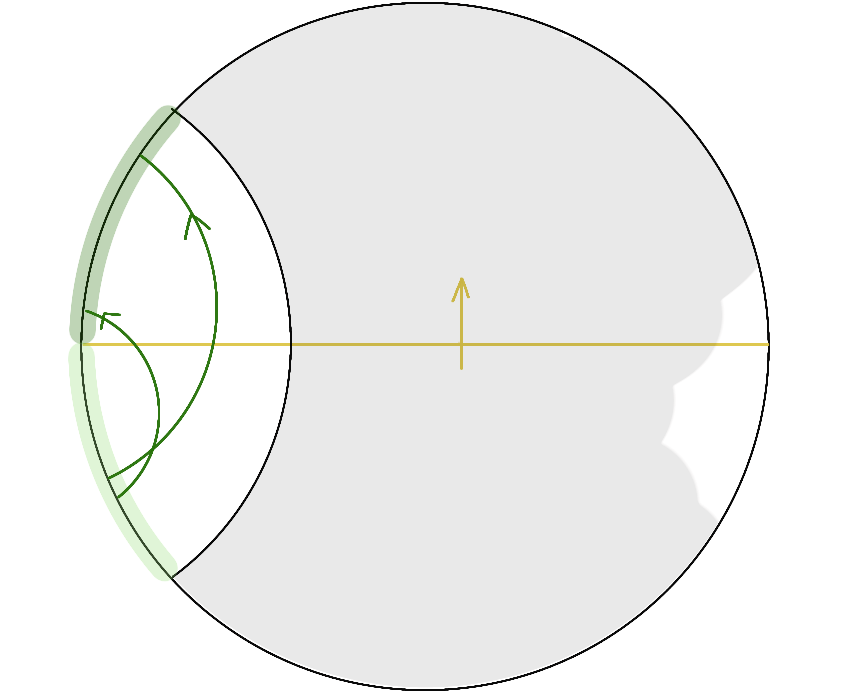}
\put(80,41){\small $a_i^j$}
\put(25,35){\small $p_i^j$}
\put(12,34){\small $p_i^{j'}$}
\put(25,50){\tiny $\eta_i^j$}
\put(13,45){\tiny $\eta_i^{j'}$}
\end{overpic}
\caption{All choices of $p_i \in a_i$ and $\theta_i$ such that (the lift of) the geodesic determined by $p_i$ and $\theta_i$ lies outside the convex core determine a proper cocycle.}\label{open}
\end{figure}

For each choice of $\underline p =(p_i(t_i))_i$ and $\underline \theta =(\theta_i(t_i, s_i))_i,$ the strip system $(\underline a, \underline p, \underline \theta) $ determines a proper affine action by Proposition \ref{good}. Consider the positive linear span $\mathcal L$ of all of these cocycles. Each cocycle in the positive linear span again determines a proper affine action, as a positive combination of cocycles with uniformly positive Margulis invariant will still have uniformly positive Margulis invariant.

Because we can independently realize an open cone's worth of translations for each separate generator with elements in $\mathcal L$, the cone $\mathcal L$ contains an open cone of $\sigma_{4n-1}(\Gamma)$-cocycles.
\endproof
\subsection{Connections to previous results}\label{prev}
\begin{figure}[h]
\centering
\begin{overpic}[scale=0.15,percent]{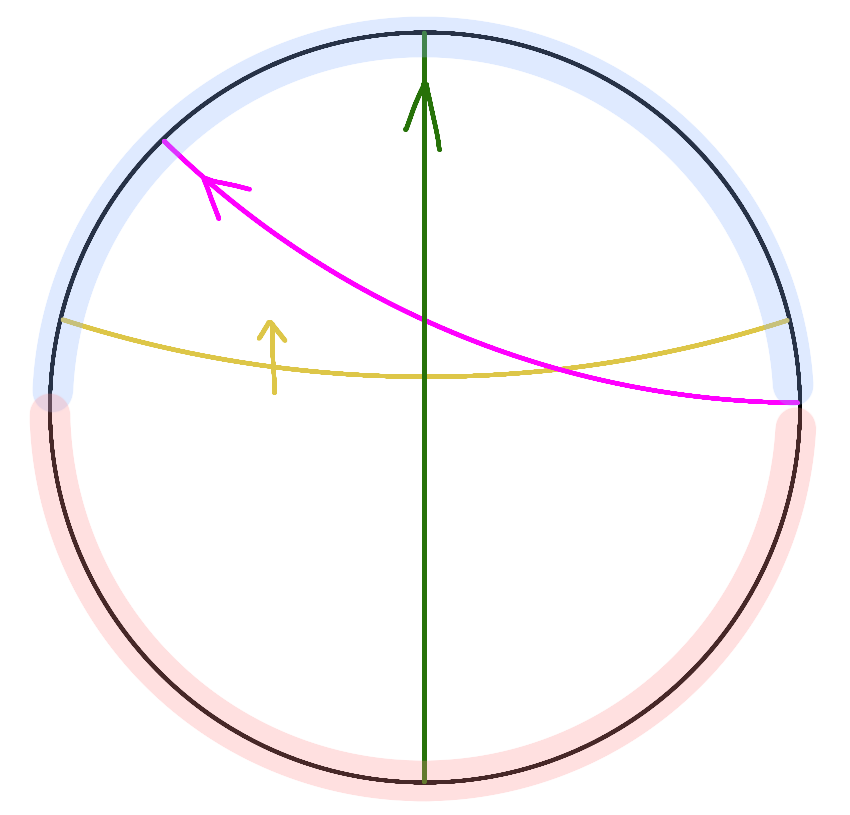}
\put(-7,46){\small $-1$}
\put(96, 46){\small $1$}
\put(80,43){\small $\gamma$}
\put(15,80){\small $t$}
\put(20,50){\small $a_i$}
\put(48,-5){\small $0$}
\put(47,94){\small $\infty$}\end{overpic}
\caption{For $n = 1,$ when we choose translations perpendicular to an arc $a_i,$ all possible contributions to the Margulis invariant are positive, as there is no sign-switching behavior  in that situation.  Were we to choose $\theta \ne 0,$ we could obtain sign-switching behavior, depending on the location of the "waist" $p_i$.}\label{nje1}
\end{figure}
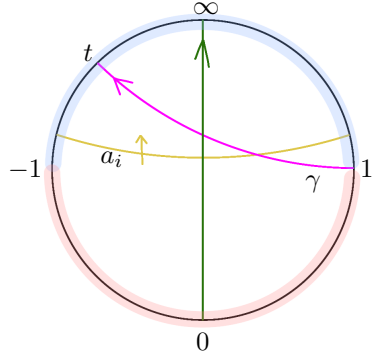

For $n=1,$ our notion of higher strip deformation coincides with \emph{infinitesimal strip deformations} from \cite{arc}. Their proof directly constructs a fundamental domain for the affine action from the strip data. It also shows much more; in the case of $n=1$ they obtain a complete parametrization of the cone of proper actions by fixing $p_i$. In that case we can choose the points $p_i$ arbitrarily as long as we choose $\theta_i = 0$ for all $a_i \in \underline a$.  We can see this from our construction by observing that the function $\rf_1(t) = \frac{t+1}{t-1}$  is positive for $t < -1$ and $t > 1,$ and negative otherwise. If we pick $\theta_i = 0,$ every curve $\gamma$ crossing $a_i$ does so with the cross-ratio of the endpoints of $\gamma$ and $\eta_i$  in $(-\infty, -1)\cup (1, \infty]$, as demonstrated in Figure \ref{nje1}. That is precisely the set where $\rf_1$ is positive.

For higher $n$, we observe sign-switching behavior by placing waists inside the convex core, even when choosing all $\theta_i$ to be $0$. But with some control over the geometry of the surface $S$ and the choice of arcs $a_i$, we have more freedom in choosing the locations of the waists $p_i.$

As in the proof of Theorem \ref{open}, suppose we choose arcs $a_i$ cutting $S$ into a topological disk, and generators $\gamma_i$ of $\Gamma$ such that each geodesic representative of $\gamma_i$ intersects only $a_i$ and does so exactly once, with the transverse orientation on $a_i$ matching the orientation of $\gamma_i.$ If we choose $p_i = \gamma_i \cap a_i$ and $\theta_i$ to be the angle between $\gamma_i$ and the perpendicular to $a_i,$ we recover Smilga's deformations from \cite{s}, which we will call of \emph{Smilga type}: 
\begin{defn}An affine deformation of a free group $\Gamma = \langle \gamma_1, \ldots, \gamma_k\rangle$ is of \emph{Smilga type} if for each generator $\gamma_i,$ its assigned translation is its neutral vector; that it, $u(\gamma_i) = x^0(\sigma_{4n-1}(\gamma_i)).$\end{defn} Smilga shows that if a dynamic expansion condition on all the $\gamma_i$ is satisfied, a deformation of Smilga type is proper. We give a similar condition on the geometry of the convex core of $S$ instead of the dynamic properties of the generators of $\Gamma.$ Further, we show that some kind of condition on the geometry or expansion really is required, by finding a deformation of Smilga type failing these conditions and also failing properness.

\begin{figure}[h!]
\centering
\begin{minipage}{0.48\textwidth}
\begin{overpic}[width=0.85\textwidth,percent]{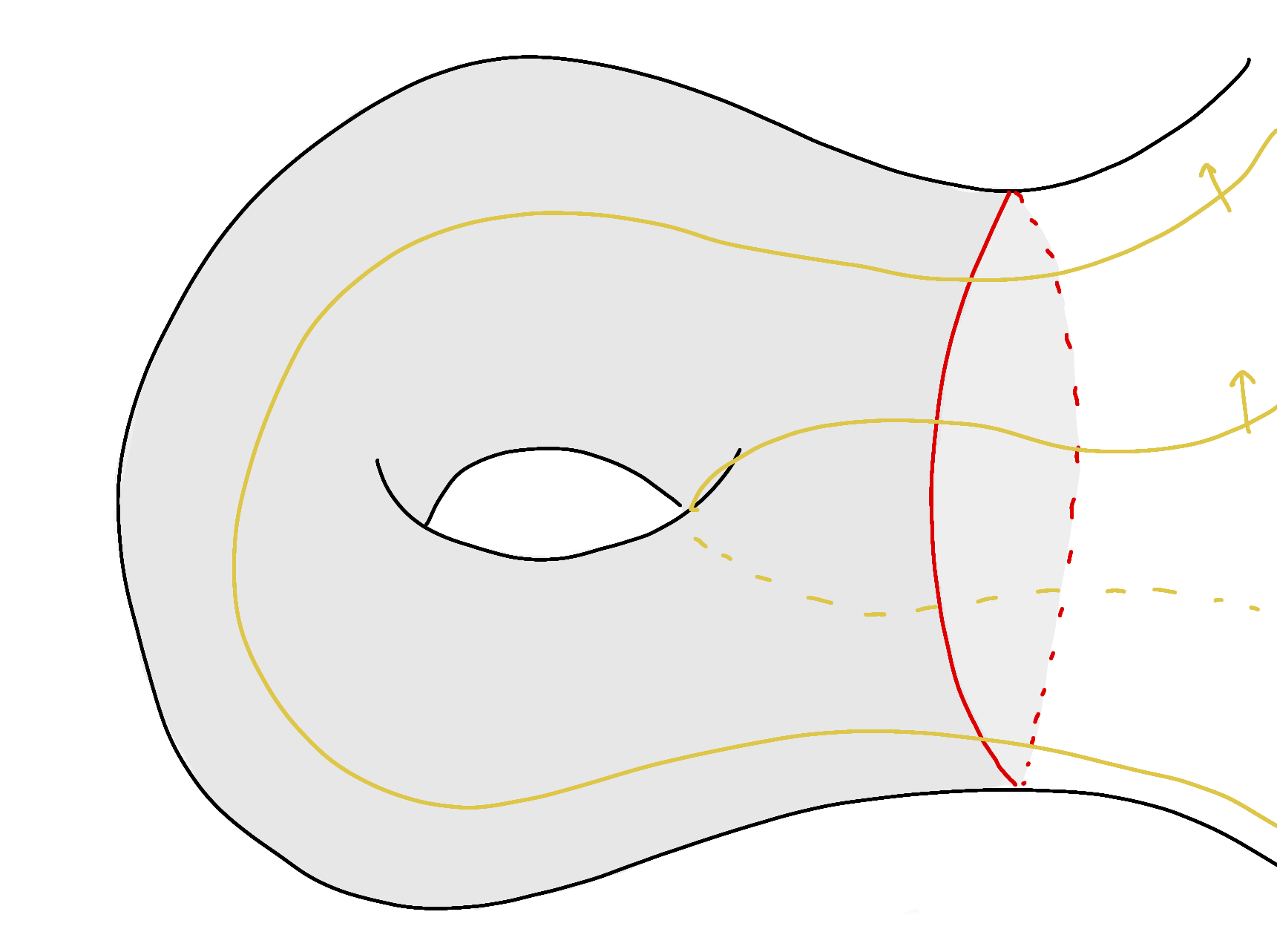}
\put(40,25){\small $\Core(S)$}
\put(75,28){ \small $\beta = \partial \Core(S)$}
\put(15,30){\small $a_1$}
\put(90,41){\small $a_2$}\end{overpic}
\end{minipage}
\begin{minipage}{0.48\textwidth}
\begin{overpic}[width=0.85\textwidth,percent]{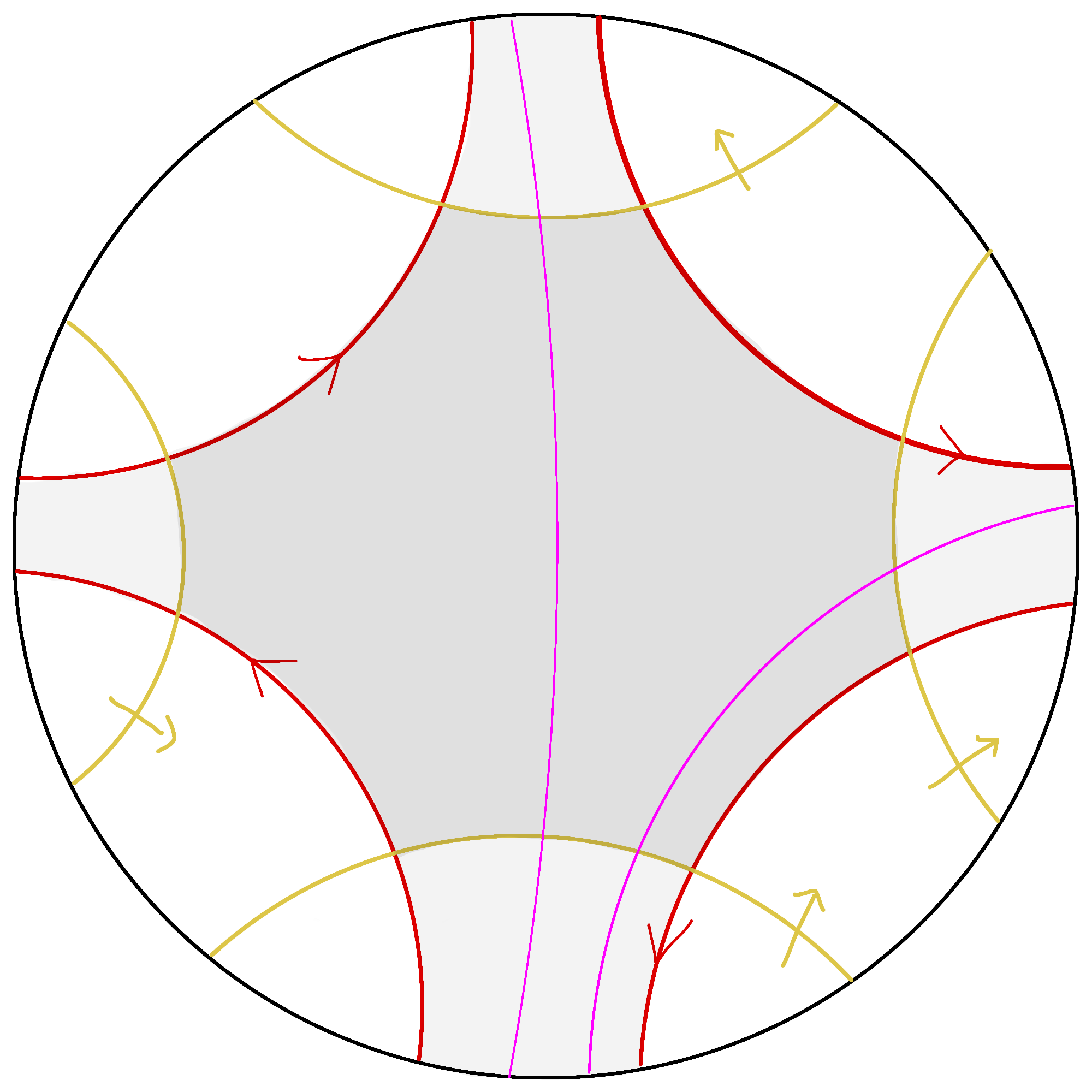}
\put(35,50){\small $K$}
\put(70,70){\small $\beta_1$}
\put(70,25){\small $\beta_2$}
\put(25,30){\small $\beta_3$}
\put(24,68){\small $\beta_4$}
\put(75,13){\small $\tilde a_2^-$}
\put(71,82){\small $\tilde a_2^+$}
\put(83,71){\small $\tilde a_1^+$}
\put(10,68){\small $\tilde a_1^-$}
\put(100,55){\small $\beta_1^+$}
\put(57,-3){\small $\beta_2^+$}
\put(-5,45){\small $\beta_3^+$}
\put(41,100){\small $\beta_4^+$}
\put(53,100){\small $\beta_1^-$}
\put(100,43){\small $\beta_2^-$}
\put(35,-3){\small $\beta_3^-$}
\put(-5,55){\small $\beta_4^-$}\end{overpic}
\end{minipage}
\caption{An example of a thin one-holed torus and its cover. The dark grey shaded region in the right picture is a fundamental domain $K$ for $\Core(S),$ and the shaded region is the lift of the convex core. All curves that stay in the convex core cross the arcs $a_i$ at an angle close to $\frac\pi2$.} \label{thint}
\end{figure}
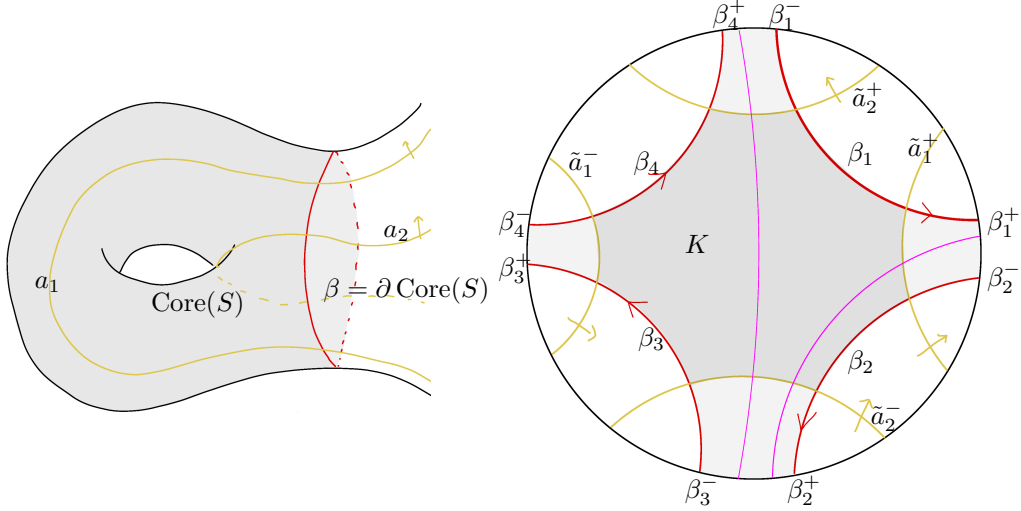

Let $K$ be a fundamental domain bounded by geodesics for $\Gamma$ acting on $\widetilde{\Core(S)} \subset \H^2$ with geodesics  $\beta_i$ being the lifts of $\partial \Core(S)$ and subarcs of the $\beta_i$ bounding $K$. Suppose the $\beta_i$ are cyclically ordered, and denote their endpoints in $\partial \H^2$ by $\beta_i^-, \beta_i^+,$ oriented in such a way that $\beta_i^- < \beta_i^+ < \beta_{i+1}^-,$ see Figure \ref{thint} for a picture. 
\begin{defn}The convex core $\Core{S}$ of $S = \Gamma \backslash \H^2$ is \emph{$c$-thin} if for all $i$, $$[\beta_{i+1}^-, \beta_{i+1}^+; \beta_i^-, \beta_i^+] > c.$$
\end{defn}
Note that  $\Core(S)$ is thinner if $c$ is bigger. A surface with large thinness is in a sense "more" convex cocompact than a surface that is less thin. The thinness condition will allow us to control the angles at which closed geodesics can intersect each other. Note that the notion of $c$-thinness is well-defined because we choose $\beta_i$ all bounding the same fundamental domain, and because the cross-ratio is $\PSL_2\R$-invariant.
Assume $a_i$ are properly embedded geodesic arcs on $S$, cutting $S$ into a topological disk, whose lifts intersect $\beta_i$ and $\beta_{i+1}.$ Denote by $\tilde a_i^{\pm}$ the lifts of the arcs $a_i$ bounding the fundamental domain $K$.  Such situations arise in the discussion above and in the proof of Theorem \ref{open}. For each $a_i,$ choose a waist inside the core of $S$ and an angle $\theta_i$ such that the axis of $\eta_i$ stays in $\Core(S).$

\begin{figure}[h]
\centering
\begin{overpic}[scale = 0.17,percent]{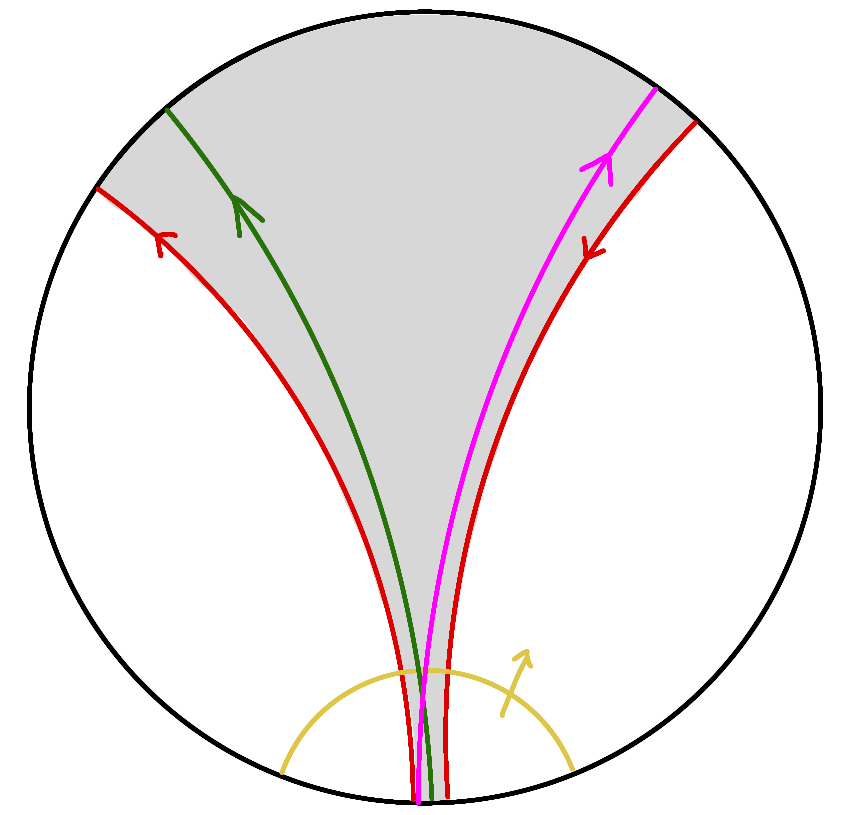}
\put(62,70){\small $\gamma$}
\put(31,71){\small $\eta_i$}
\put(60,40){\small $\beta_{i}$}
\put(25,40){\small $\beta_{i+1}$}
\put(83, 82){\small $\beta^-_i$}
\put(54,-5){\small $\beta^+_i$}
\put(40,-5){\small $\beta^-_{i+1}$}
\put(1,75){\small $\beta^+_{i+1}$}
\put(63,13){\small $\tilde a_i$}\end{overpic}
\caption{If two consecutive arcs $\beta_i, \, \beta_{i+1}$ are close, that restricts the possible angles between curves lying in the shaded "wedge" they make.}
\label{wedge}
\end{figure}
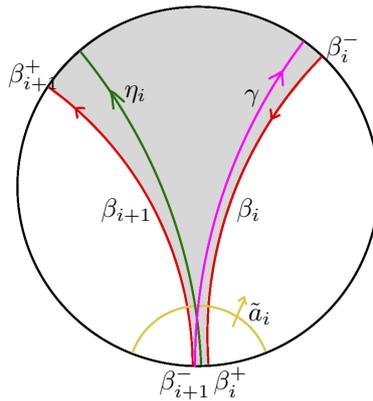

\begin{prop}\label{thin}For every $n$ there exists a $c_n$ such that if the convex core of $\Gamma \backslash \H^2$ is at least $c_n$-thin, the affine deformation determined by the strip data $(\underline a, \underline p, \underline \theta)$ as above gives a proper affine deformation of $\sigma_{4n-1}(\Gamma).$
\end{prop}
\proof Let $x_n$ be the zero of largest modulus of $\rf_n(t)$ and let $c_n > |x_n|.$ Let $\eta_i$ be a translation whose axis intersects $a_i$ and stays within $\Core(S),$ agreeing with the transverse orientation on $a_i.$ One of the fixed points of $\eta_i$ lies between $\beta_i^+$ and $\beta_{i+1}^-,$ and the other one lies between $\beta_i^-$ and $\beta_{i+1}^+.$ The same is true for any $\gamma \in \Gamma$ intersecting $a_i$. The largest negative value a cross-ratio between the endpoints of such a pair of curves can take is obtained by the curve connecting $\beta_i^+$ to $\beta_{i+1}^-$ and the curve between $\beta_i^-$ and $\beta_{i+1}^+$, see Figure \ref{wedge}. That cross-ratio is $-[\beta_{i+1}^-, \beta_{i+1}^+; \beta_i^-, \beta_i^+]$, which is smaller than $-c_n$ because of the thinness condition. Therefore, by Proposition \ref{invt1}, the Margulis invariant contribution upon each crossing of an arc in $\underline a$ by $\gamma$ is at least $\rf_n(c_n) > 0$.

Similarly to the argument in the proof of Theorem \ref{th}, we get uniformity of the Margulis invariant from the fact that $S$ is convex cocompact, and thus each closed curve crosses an arc in $\underline a$ least once per some fixed length, gaining a uniform contribution of at least $\rf_n(c_n)$ upon each crossing. Then by \cite{GLM}, we can conclude that this action is proper.
\endproof
\begin{remark}The thinness parameter $c_n$ gets bigger as $n$ grows, so we need more and more control over the geometry of the surface - or the dynamical properties of the generators of $\Gamma$ - for the Smilga construction to work. On the other hand, placing the waists and axes of $\eta_i$ outside of the convex core provides a strip system determining a proper cocycle in any dimension.
\end{remark}

Below, we give a concrete example of a one-holed torus group, where the convex core is not sufficiently thin and the deformation of Smilga type for $n=3$ does not determine a proper action.
\begin{ex}For $\lambda > 1 + \sqrt2$ the group $\Gamma$ generated by the matrices $$X = \begin{bmatrix}\lambda&0\\0&\frac{1}{\lambda}\end{bmatrix}, \,Y=\frac12\begin{bmatrix}\lambda +\frac1\lambda & \lambda - \frac1\lambda\\\lambda-\frac1\lambda & \lambda + \frac1\lambda\end{bmatrix}$$ is a discrete free group with quotient $S=\Gamma \backslash \H^2$  a convex cocompact one-holed hyperbolic torus. It is very symmetric in the sense that the lengths of its standard generators are the same.

Describe the deformation of $\sigma_{11}(\Gamma) < \SO(6,5)$ by setting $u(X) = x^0(\sigma_{11}(X))$ and $u(Y) = x^0(\sigma_{11}(Y)).$ This is a cocycle of Smilga type. As a higher strip deformation, it is obtained from the strip system $(\underline a, \underline p, \underline \theta)$ where the arcs of $\underline a$ are the two arcs cutting $S$ up into a disk, with $a_1$ perpendicular to the axis of $X$ and the axis of $Y$ perpendicular to $a_2$, with $p_i$ the points of intersection and $\theta_i = 0.$
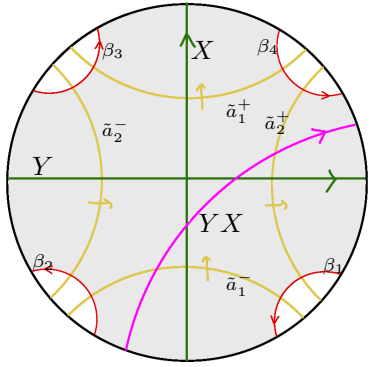
\begin{figure}[h]
\centering
\begin{overpic}[scale=0.17,percent]{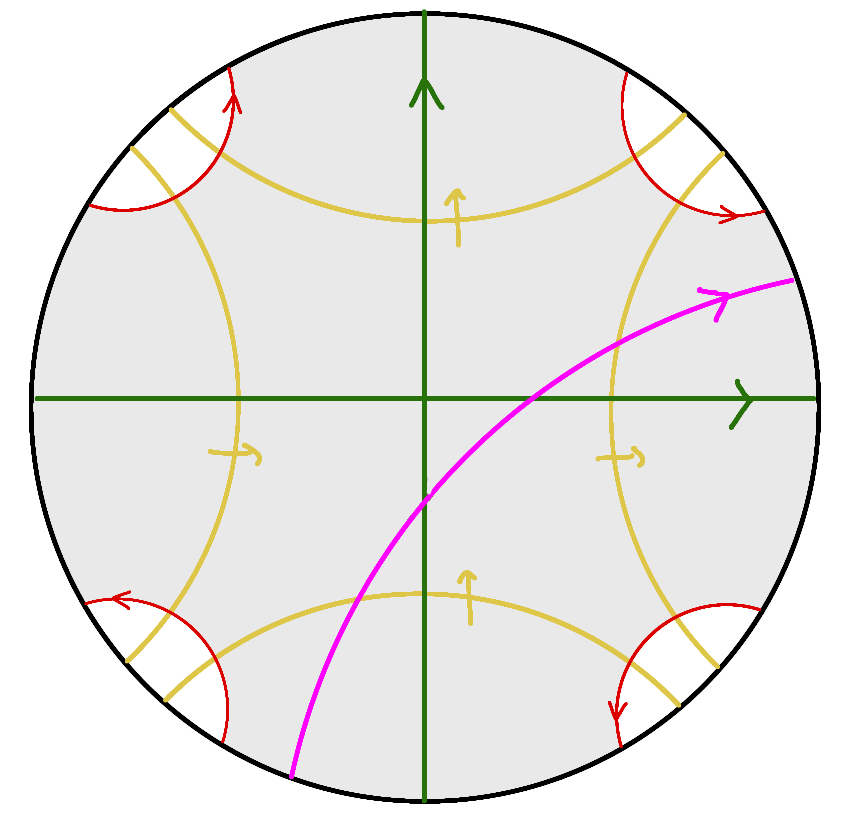}
\put(51,80){\small $X$}
\put(10,50){\small $Y$}
\put(53,35){\small $YX$}
\put(68,82){\tiny $\beta_4$}
\put(85,25){\tiny$\beta_1$}
\put(10,26){\tiny$\beta_2$}
\put(28,81){\tiny$\beta_3$}
\put(60,20){\tiny$\tilde a_1^-$}
\put(60,65){\tiny$\tilde a_1^+$}
\put(28,60){\tiny$\tilde a_2^-$}
\put(70,62){\tiny$\tilde a_2^+$}\end{overpic}
\caption{A very symmetric one-holed torus used in this example.}\label{torus}
\end{figure}

No matter what $\lambda$ is, the Margulis invariants of $X$ and $Y$ are $1$; $\alpha_u(X) = B(x^0(X) , u(X)) = B(x^0(X), x^0(X)) = 1$ and the same for $Y$. However, the Margulis invariant of $YX$ depends on $\lambda.$

The geodesic representative of $YX$ crosses $a_1$ once and $a_2$ once, both at the same angle. Due to this symmetry, it is enough to compute the Margulis invariant contribution to $YX$ as it crosses $a_1$ once. The cross-ratio of the endpoints of $X$ and the endpoints of $YX$ is
$$[X^+, X^-;(YX)^+, (YX)^-] = [\infty,0; (YX)^+, (YX)^-] = \frac{1 + \lambda + \sqrt{1 + 6\lambda^2 + \lambda^4}}{1 + \lambda - \sqrt{1 + 6\lambda^2 + \lambda^4}}.$$

Plugging the expression into $\rf_3(t),$ we get 
$$\frac{(1 + \lambda^2) (1 - 16 \lambda^2 - 4 \lambda^4 - 16 \lambda^6 + \lambda^8)}{(1 + 6 \lambda^2 + \lambda^4)^{\frac52}},$$
whose graph for $\lambda > 1+\sqrt2$ is
\begin{figure}[h!]
    \centering
        \begin{overpic}[scale=0.4,percent]{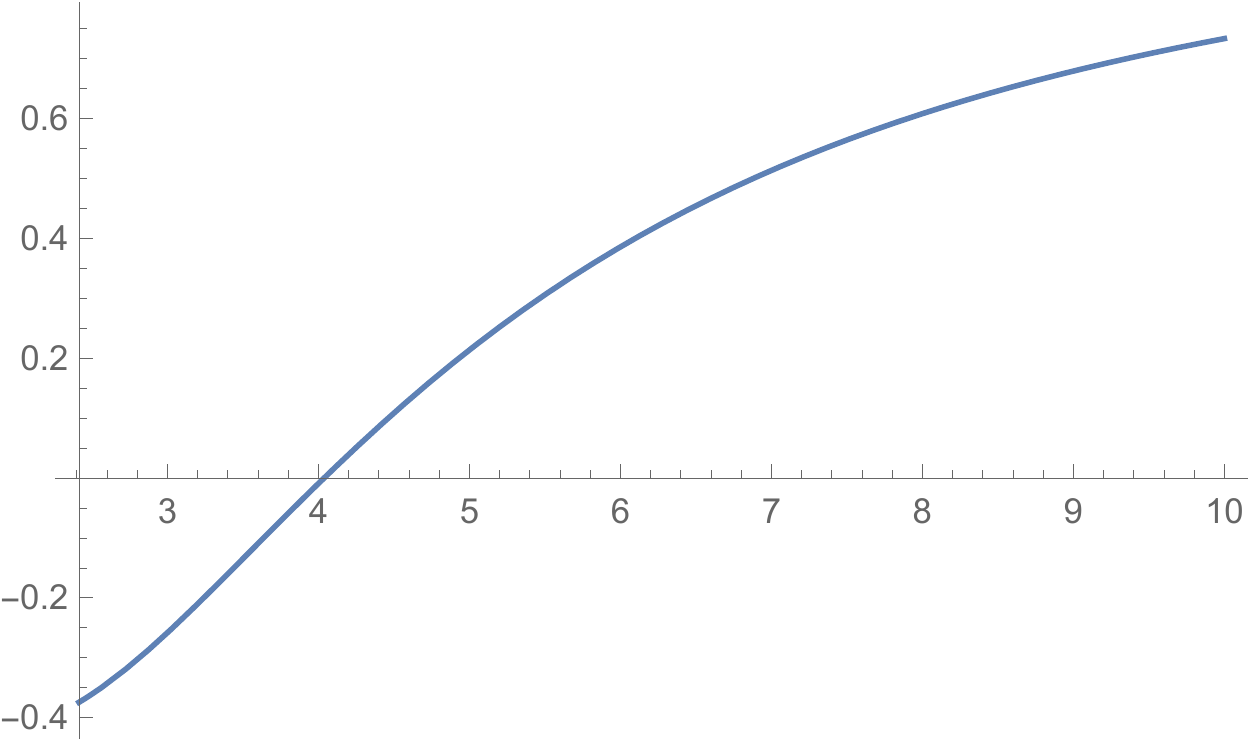}
\put(100,15){\tiny $\lambda$}
\put(-18,55){\tiny $\alpha_u(YX)$} \end{overpic}
   \caption{The graph of the Margulis invariant of $YX$ in $\R^{6,5}$ upon one crossing with respect to $\lambda$.}
\end{figure}

We see that for say, $\lambda = 3,$ the Margulis invariant of $YX$ is negative, but as the Margulis invariant of $X$ is positive, this action can't be proper.  This example shows that the thinness condition (or Smilga's expansion condition on the generators) is not in vain, and again contrasts  with the situation for $n = 1,$ where the strip data determining the deformation gives a proper action for all $\lambda > 1 + \sqrt2$.

 Even when $\lambda$ is large enough to guarantee positivity of the Margulis invariant of $YX,$ that is not necessarily enough to ensure properness of the action of $\Gamma$, as some other element could still have negative Margulis invariant. Increasing $\lambda$ enough to ensure thinness of the convex core as per Proposition \ref{thin} would give us positive Margulis invariant for all curves at once. In this case, we can directly compute the thinness of $\Gamma \backslash \H^2,$ as the boundary of $\Core(S)$ is $$XYX^{-1}Y^{-1}=\begin{bmatrix}1 + \lambda^2 + 3\lambda^4 -\lambda^6 & (\lambda^2 -1)^2(1+\lambda^2) \\ -(\lambda^2-1)^2(1+\lambda^2)& -1+3\lambda^2 + \lambda^4 + \lambda^6)\end{bmatrix},$$ so its axis and the axes of its conjugates in $\H^2$ are the geodesics $\beta_i$ from the thinness condition. Because of the symmetry of the once-punctured torus in this example, it is enough to compute the thinness with respect to two consecutive arcs $\beta_2, \beta_1.$ Their endpoints respectively are $\frac12 (-1 + \lambda^2 + \sqrt{1-6\lambda^2 + \lambda^4}), \, \frac12 (-1 + \lambda^2 - \sqrt{1-6\lambda^2 + \lambda^4})$ and  $ \frac{2}{-1 + \lambda^2 - \sqrt{1-6\lambda^2 + \lambda^4}}, \, \frac{2}{-1 + \lambda^2 + \sqrt{1-6\lambda^2 + \lambda^4}}.$ Their cross-ratio is $$\frac{(\lambda^2 -1)^2}{4\lambda^2}.$$ The smallest zero of $\rf_3$ is $t_0= \frac12(-12 -\sqrt{70} - \sqrt{6(35 + 4 \sqrt{70})}),$ which is approximately $-20.3174.$  Proposition \ref{thin} ensures that if we choose $\lambda$ such that $\frac{(\lambda^2 -1)^2}{4\lambda^2}> t_0,$ the affine deformation determined by the strip data will be proper. The inequality is first satisfied for $\lambda  \approx 9.12456.$ 
\end{ex}

\section{Actions of virtually free groups}\label{virt}

We can use the Margulis invariant and our construction to make affine deformations of virtually free groups. As proof of concept, we present actions of the free product of the cyclic groups of order $2$ and $3$ on $\R^7.$ Throughout this section, denote by $\Gamma $ a subgroup of $\PSL_2\R$ abstractly isomorphic to $C_2 \star C_3.$ We can generate $\Gamma$ with   a degree 2 rotation $S$ and a a degree 3 rotation $R$ in the hyperbolic plane. Up to conjugation, such an embedding of $C_2 \star C_3$ is determined by the distance between the fixed points of the two generating rotations. It is an embedding if this distance is at least $\ln\sqrt 3,$ the radius of an inscribed circle of an ideal triangle. The quotient $\Gamma \backslash \H^2$ is an orbifold with two cone points, one of degree two and one of degree three. It is covered by a thrice-punctured sphere, with the fundamental group $\Gamma'$ of the thrice-punctured sphere generated by $(SR)^2$ and $(RS)^2.$ If we choose the rotations $S$ and $R$ to be $\begin{bmatrix} 0 & 1\\-1&0 \end{bmatrix}$ and $\begin{bmatrix}0&-1\\ 1 & 1 \end{bmatrix},$ we get the standard $\PSL_2\Z < \PSL_2\R.$

If we want to describe a cocycle of $\sigma_{4n-1}(\Gamma),$ it is enough to describe a cocycle of the fundamental group of the thrice-punctured sphere covering $\Gamma\backslash \H^2$ invariant under the action of $\Sym_3=\Gamma/\Gamma'$ on the thrice-punctured sphere. In that case, the action of the free group $\Gamma'$ will extend to the action of $\Gamma,$ and because $\Gamma'$ has finite index in $\Gamma,$ properness of the $\Gamma'$ action suffices for properness of the $\Gamma$ action.  There is only one strip system on the thrice-punctured sphere that descends to a strip system  on $\Gamma\backslash \H^2.$ Its arcs are the lifts of the arc through the degree $2$ cone point on $\Gamma \backslash \H^2$, which has three preimages in the thrice-punctured sphere.  The lifts of any other properly embedded arc on $\Gamma \backslash \H^2$ would intersect in $\Gamma' \backslash \H^2.$ To preserve symmetry, the points in $\underline p$ have to be lifts of the degree 2 cone point, and the angles $\underline \theta$ all have to be $0$. We choose transverse orientations on the arcs arbitrarily; as long as all $\theta_i$ are $0$, we get an equivariant cocycle out of the strip system.

The standard embedding $\PSL_2\Z < \PSL_2\R$ is a group with parabolic elements, the affine deformations of which we cannot say much about using the methods presented here. There is work extending the notion of the Margulis invariant to groups with parabolic elements, see for instance \cite{par}. However, to the best of the author's knowledge, there is no known sufficient condition on the Margulis invariant that would guarantee properness of the action in the presence of parabolic elements. Therefore, when not stated otherwise, assume $\Gamma < \PSL_2\R$ (and thus $\Gamma'$) is convex cocompact, meaning that the distance between the fixed points in $\H^2$ of the generating rotations $R$ and $S$ is greater than $\ln 3. $

\begin{figure}[h!]
    \centering
    \begin{minipage}{0.5\textwidth}
        \centering
        \includegraphics[width=0.9\textwidth]{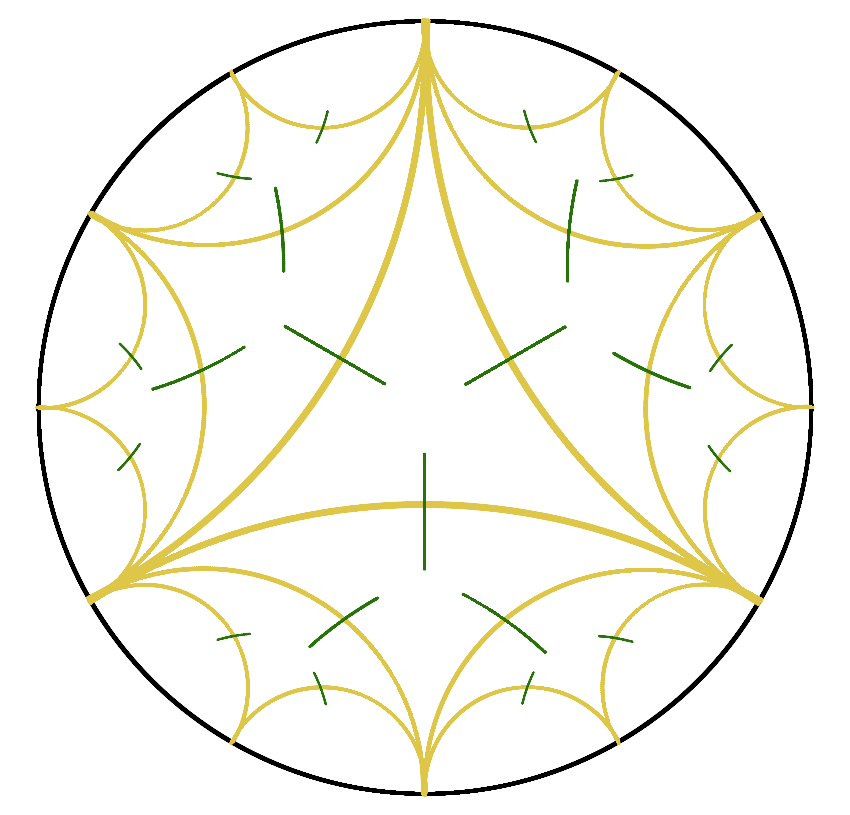} 
    \end{minipage}\hfill
    \begin{minipage}{0.5\textwidth}
        \centering
        \includegraphics[width=0.9\textwidth]{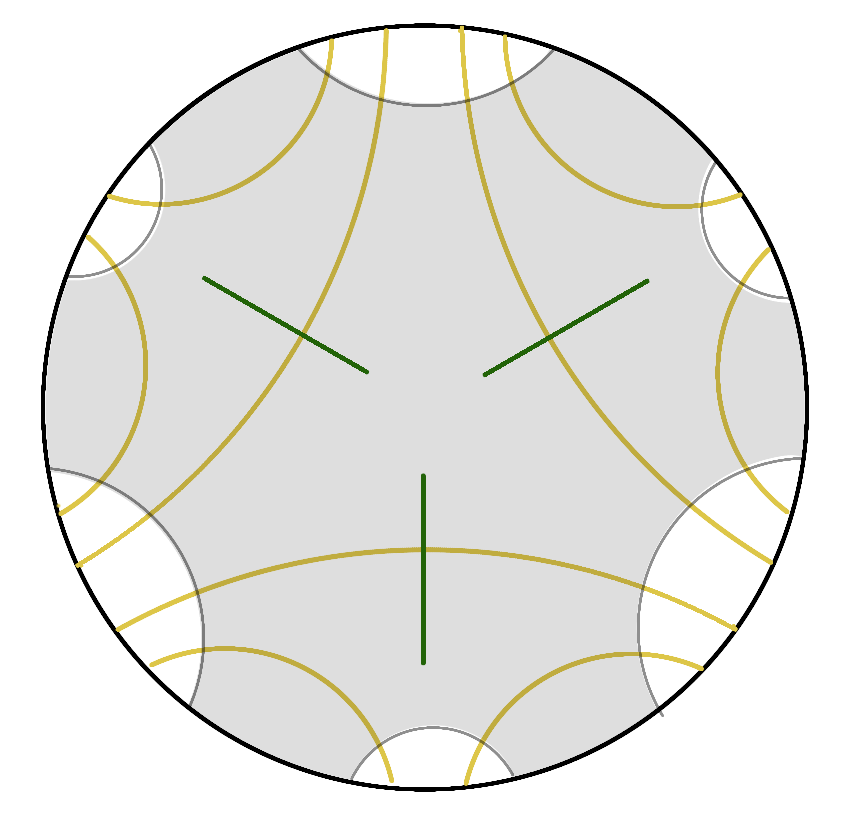} 
    \end{minipage}

\caption{The covers of the orbifold $\Gamma\backslash \H^2$, with the lift of the strip system. On the left, we have the standard way $\PSL_2\Z$ lies inside $\PSL_2\R,$ and in the right picture, we "opened up the cusp" of the modular surface.}\label{psl}
\end{figure}

We focus here on dimension 7, because the cohomology group $\operatorname H^1(\Gamma, \R^{7}) = \R$. Up to conjugacy and scaling, there is only one affine cocycle, and we can realize it with a higher strip deformation. In particular, analyzing that higher strip deformation will tell us that all non-nullcohomologous cocycles of $\Gamma$ determine a proper affine action on $\R^7.$ In higher dimensions,  the first group cohomology group is no longer one-dimensional. Analyzing  the cocycle determined by the described strip system is still possible, but does not yield complete information about all possible affine deformations. 

Let us now analyze the one higher strip deformation  we have. Denote by $a_1, \, a_2, \, a_3$ the arcs on the thrice-punctured sphere connecting two distinct cusps (or funnels), and with $p_1, \, p_2,\, p_3$ the points on them that are lifts of the degree 2 cone point on $\PSL_2\Z \backslash \H^2.$ let $\theta_i = 0.$  As noted in the discussion in Section \ref{prev}, we can get negative contributions to the Margulis invariant from these strips. They are neither of Smilga type - the translations $\eta_i$ do not correspond to any closed loops on $\Gamma' \backslash \H^2$ -, nor do the geodesics through $p_i$ lie entirely outside the convex core. However, in dimension 7, a striking coincidence happens. Every closed curve on the three-holed sphere will have to cross an even number of arcs in $\underline a,$ and we will show that the contribution to the Margulis invariant of a cocycle with values in $\R^{4,3}$ upon two consecutive crossings is always positive!

\begin{prop}\label{angles}Let $  C_2\star C_3\cong \Gamma < \PSL_2\R$ be convex cocompact. Let $\Gamma' < \Gamma$ be the index 6 free subgroup and let  $S = \Gamma' \backslash \H^2$ be the three-holed sphere. Let the axis of a loxodromic element $\gamma \in \Gamma'$ cross the arcs $\tilde a_1$ and $\tilde a_2$ in that order. For $ i = 1,2$, let $\eta_i$ be the loxodromic translation with axis through $\tilde p_i \in \tilde a_i$ with transverse angle $\theta = 0.$ Then $B(x^0(\sigma_7(\gamma)), x^0(\sigma_7(\eta_1)) + x^0(\sigma_7(\eta_2))) > 0.$\end{prop}
\proof

Recall from Proposition \ref{invt1} that for a higher strip deformation in $\R^{4,3},$ the contribution to the Margulis invariant when $\gamma$ crosses an arc $a_i$ once is given in terms of the the cross-ratio $t$ between the endpoints of $\gamma$ and $\eta_i$ by $\rf_2(t) = \frac{1}{(t-1)^3}\left(1 + 9 t + 9t^2 + t^3 \right).$  Because the orientations of $\gamma$ and $\eta_i$ will always agree, we only need to worry about what happens for $t \in (-\infty, 1 ) \cup (1, \infty).$ On that set, $\rf_2$ attains a minimum at 
$\frac{-3-\sqrt 5}{2},$ and 
$\rf_2(\frac{-3-\sqrt5}{2}) = -\frac{1}{\sqrt 5}.$

If the axis of $\gamma$ doesn't cross the axes of $\eta_1$ or $\eta_2,$ we know from Lemma \ref{bigt} that $B(x^0(\sigma_7(\gamma)), x^0(\sigma_7(\eta_i)) > 1,$ so $B(x^0(\sigma_7(\gamma)), x^0(\sigma_7(\eta_1)) + x^0(\sigma_7(\eta_2))) > 2 > 0.$

If the axis of $\gamma$ intersects the axis of $\eta_1$ but not of $\eta_2,$ then $B(x^0(\sigma_7(\gamma)), x^0(\sigma_7(\eta_1)) + x^0(\sigma_7(\eta_2))) > -\frac1{\sqrt5} + 1 > 0.$

The remaining case is when the axis of $\gamma$ intersects the axes of both $\eta_1$ and $\eta_2$.
Then, these three axes determine a triangle in the hyperbolic plane. Denote by $\phi_1$ the angle between the axes of $\gamma$ and $\eta_1,$ and by $\phi_2$ the angle between the axes of $\gamma$ and $\eta_2$. 
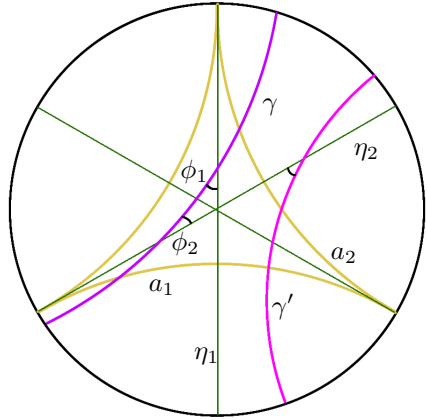
\begin{figure}[h!]
\centering
\begin{overpic}[scale=0.2,percent]{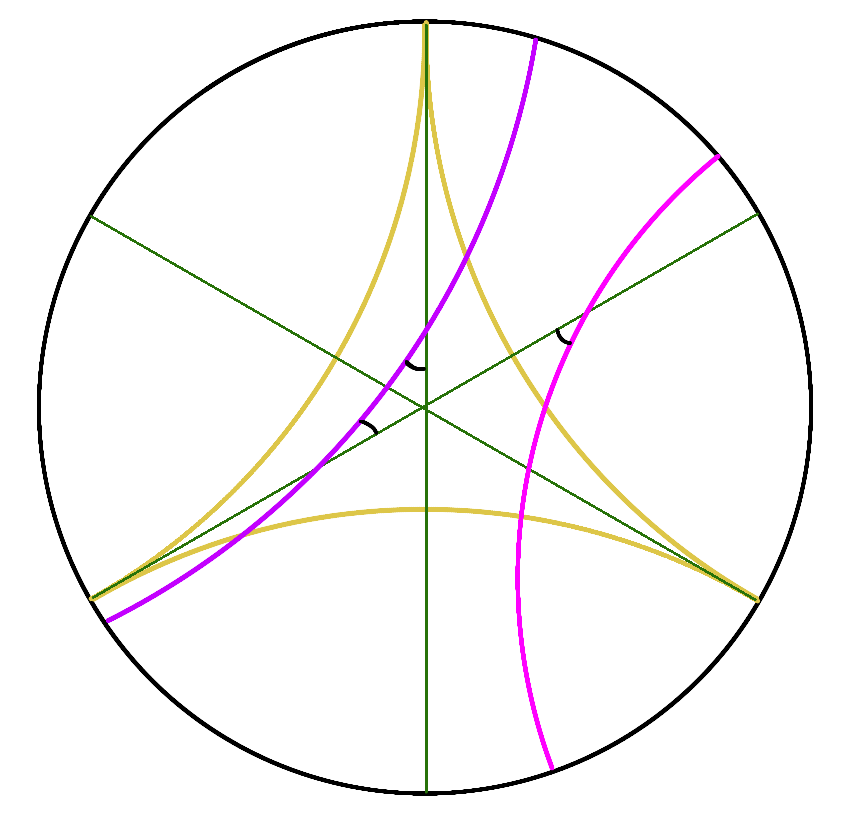}
\put(62,25){\small $\gamma'$}
\put(60,70){\small $\gamma$}
\put(35,30){\small $a_1$}
\put(45,15){\small $\eta_1$}
\put(80,60){\small $\eta_2$}
\put(75,37){\small $a_2$}
\put(42,55){\small $\phi_1$}
\put(40,39){\small $\phi_2$}
\end{overpic}
\caption{Two of the ways an element $\gamma$ can cross two consecutive arcs in the cover of $S$.}\label{cross}
\end{figure}

\begin{figure}[h]
\includegraphics[scale=0.6]{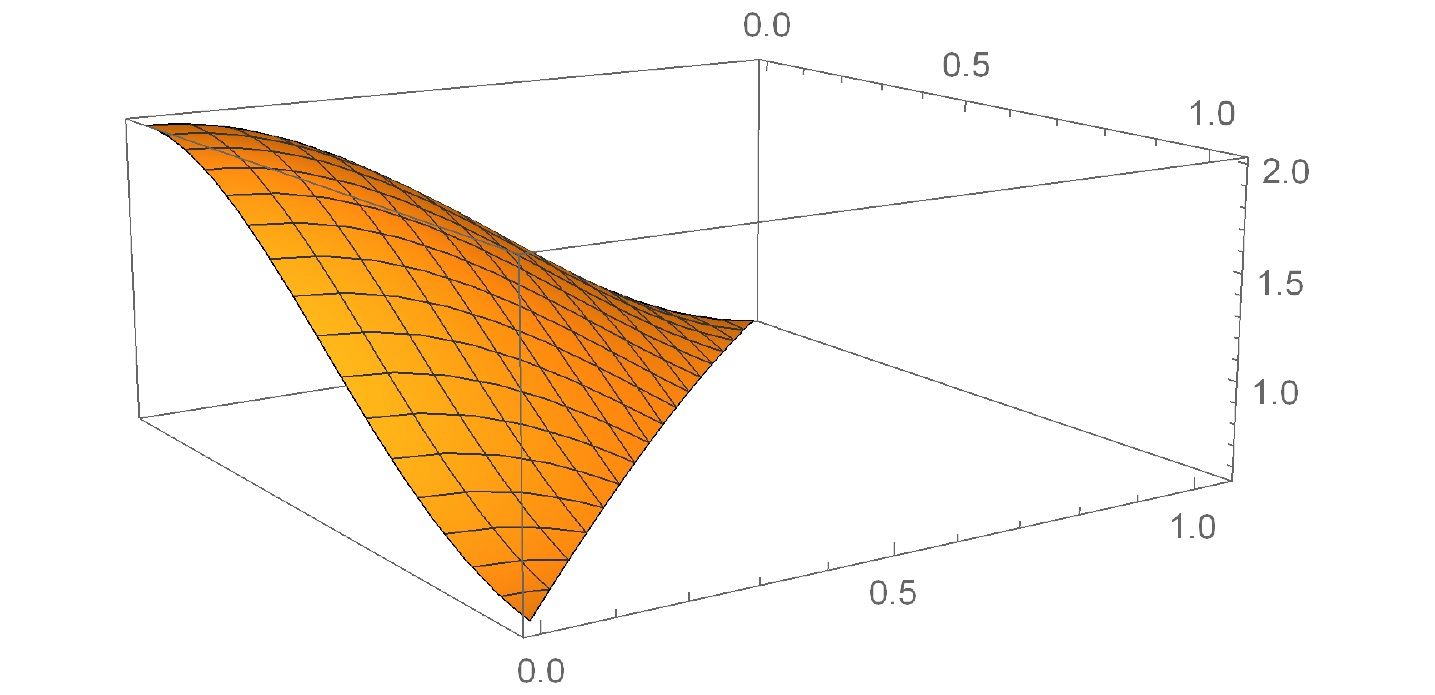}
\caption{The graph of $s$.}\label{s}
\end{figure}

Recall from Remark \ref{withangles} that when the axes intersect at angle $\phi$, we can rewrite $\rf_2$ in terms of $\phi$ as $\rf_2(\phi) = \frac18(3 \cos \phi + 5 \cos 3\phi).$ Then we can write the Margulis invariant contribution coming from $\gamma$ crossing $\tilde a_1$ and $\tilde a_2$ as $s(\phi_1, \phi_2) = \frac18(3 \cos \phi_1 + 5 \cos 3\phi_1)+\frac18(3 \cos \phi_2 + 5 \cos 3\phi_2).$ The graph of the function $s$ is sketched in Figure \ref{s}.
\endproof
Due to the symmetry of the arc system, the angle between the axes of $\eta_1$ and $\eta_2$ is $\frac{2\pi}{3}$. As the sum of angles in a triangle in the hyperbolic plane is less than $\pi,$ it follows that $\phi_1 + \phi_2 < \frac\pi3.$ It is a straightforward calculus problem to show that on the domain determined by $0 \le \phi_1, \phi_2$ and $\phi_1 + \phi_2 \le \frac \pi3,$ the function $s$ is positive. It attains its minimum of  $\frac{9}{16}$ only when one of the angles $\phi_1, \phi_2$ is $0$ and the other is $\frac \pi3,$ which can only happen if the axis of one of the $\eta_i$ coincided with that of $\gamma.$

By \cite{GLM}, the preceding proof shows:

\begin{corollary} Let $\Gamma \cong C_2\star C_3$ be a convex cocompact deformation of $\PSL_2\Z$ in $\PSL_2\R.$ Then every non-trivial affine deformation of $\sigma_7(\Gamma)$ determines a proper affine action on $\R^{7}.$
\end{corollary}

\bibliographystyle{amsalpha}
\bibliography{bib}

\end{document}